\documentclass[12pt]{amsart}
\usepackage{amscd,amssymb,verbatim}
\usepackage{amssymb,amsmath,amsthm,epsfig}

\newcommand{\N}{\mathbb{N}}

\newcommand{\R}{\mathbb{R}}
\newcommand{\supp}{\mbox{supp}}

\theoremstyle{plain}
\newtheorem{Thm}{Theorem}[section]
\newtheorem{Prop}[Thm]{Proposition}
\newtheorem{Lem}[Thm]{Lemma}
\newtheorem{Cor}[Thm]{Corollary}
\newtheorem{Rmk}[Thm]{Remark}

\newtheorem{Qst}[Thm]{Question}
\newtheorem{Def}[Thm]{Definition}

\theoremstyle{remark}

\begin{document}

\title{On the structure of the spreading models of a Banach space}
\author{G.~Androulakis, E.~Odell,
Th.~Schlumprecht, N. Tomczak-Jaegermann}
\thanks{This research was supported by NSF, NSERC and the Pacific
Institute for the Mathematical Sciences.  In addition, the fourth
author holds the Canada Research Chair in Mathematics.}

\date{}
\maketitle

\noindent
{\bf Abstract} We study some questions concerning the structure of the
set of spreading models of a separable infinite-dimensional Banach
space $X$. In particular we give an example of a reflexive $X$ so that
all spreading models of $X$ contain $\ell_1$ but none of them is
isomorphic to $\ell_1$. We also prove that for any countable set $C$
of spreading models generated by weakly null sequences there is a
spreading model generated by a weakly null sequence which dominates
each element of $C$. In certain cases this ensures that $X$ admits,
for each $\alpha < \omega_1$, a spreading model $(\tilde
x_i^\alpha)_i$ such that if $\alpha < \beta$ then $(\tilde
x_i^\alpha)_i$ is dominated by (and not equivalent to)
$(\tilde x_i^\beta)_i$. Some applications of these ideas are used to
give sufficient conditions on a Banach space for the existence of a
subspace and an operator defined on the subspace, which is not a
compact perturbation of a multiple of the inclusion map.

\section{Introduction}\label{sec1}

It is known that for every normalized basic
sequence $(y_i)$ in a Banach space and for every $\varepsilon_n \searrow 0$
there exists a subsequence $(x_i)$
and a normalized basic sequence $(\tilde{x}_i)$ such that: \ For all
$n\in \mathbb{N}$, $(a_i)^n_{i=1}\in [-1,1]^n$
and $n\le k_1 <\ldots <k_n$
\begin{equation}\label{eq0}
\left|\bigl\| \sum^n_{i=1} a_ix_{k_i}\bigr\| - \bigl\|\sum^n_{i=1}
a_i\tilde x_i\bigr\|\right| <\varepsilon_n.
\end{equation}
The sequence $(\tilde x_i)$ is called the {\em spreading model of
$(x_i)$} and it is a suppression-1 unconditional basic sequence if
$(y_i)$ is weakly null (see \cite{BS1} and \cite{BS2}; see also
\cite{BL}(I.3.Proposition 2) and \cite{O} for more about spreading 
models).  This in
conjunction with Rosenthal's $\ell_1$ theorem \cite{R}, yields that
every separable infinite dimensional Banach space $X$ admits a
suppression 1-unconditional spreading model $(\tilde x_i)$. In fact
one can always find a $1$-unconditional spreading model \cite{R2}.  It
is natural to ask if one can always say more. What types of spreading
models must always exist?
Sometimes we refer to $\widetilde X = [\tilde x_i\colon \ i \in
\mathbb{N}]$, the closed linear span of $(\tilde x_i)$, as the
spreading model of $(x_i)$. By James' well known theorem \cite{J}
every such $X$ thus admits a spreading model $\widetilde X$ which is
either reflexive or contains an isomorph of $c_0$ or $\ell_1$. It was
once speculated that for all such $X$ some spreading model $(\tilde
x_i)$ must be equivalent to the unit vector basis of $c_0$ or $\ell_p$
for some $1\le p <~\infty$ but this was proved to be false \cite{OS}.
A replacement conjecture was brought to our attention by V.D.~Milman:\
must every separable space $X$ admit a spreading model which is either
isomorphic to $c_0$ or $\ell_1$ or is reflexive? In section~2 we show
this to be false by constructing a space $X$ so that for all spreading
models $\widetilde X$ of $X$, $\widetilde X$ contains $\ell_1$ but
$\widetilde X$ is never isomorphic to $\ell_1$. The example borrows
some of the intuition behind the example of \cite{OS}. That space had
the property that amongst the $\ell_p$ and $c_0$ spaces only $\ell_1$
could be block finitely representable in any spreading model $(\tilde
x_i)$ yet no spreading model could contain $\ell_1$.

The motivation behind our example comes from the ``Schreierized''
version $S(d_{w,1})$ of the Lorentz space $d_{w,1}$. Let $1=w_1
>w_2>\ldots$ with $w_n\to 0$ and $\sum^\infty_{n=} w_n=\infty$. Then
$d_{w,1}$ is the sequence space whose norm is given by
$$\|x\| = \sum_n w_nx^*_n$$
where $x$ is the sequence $(x_n)$ and
$(x^*_n)$ is the decreasing rearrangement of $(|x_n|)$. One could then
define the sequence space $S(d_{w,1})$ as the completion of $c_{00}$
(the linear span of finitely supported sequences of reals) under
$$
\|x\| = \sup_{n\in\N,n\le k_1<k_2<\ldots<k_n } \sum^n_{i=1}
w_ix^*_{k_i}.
$$
In this case the unit vector basis $(e_i)$ has a spreading model,
namely the unit vector basis of $d_{w,1}$, which is not an $\ell_1$
basis but whose span is hereditarily $\ell_1$. $S(d_{w,1})$ is
hereditarily $c_0$ so it does not solve Milman's question. In order to
avoid $c_0$ one may also define the ``Tsirelsonized'' version
$T(d_{w,1})$ of $d_{w,1}$. $T(d_{w,1})$ is the completion of $c_{00}$
under the implicit equation
$$
\|x\| = \max\left(\|x\|_\infty, \sup \sum^n_{i=1} w_i
   \|E_ix\|^*\right)
$$
where the supremum is taken over all integers $n$, and all {\em
   admissible} sets $(E_i)^n_{i=1}$ i.e.\ $n\le E_1<\ldots< E_n$ (this
means $n\le \min E_1 \le \max E_1< \min E_2\le\ldots)$. $E_ix$ is the
restriction of $x$ to the set $E_i$. It may well be that $T(d_{w,1})$
has the properties we desire but we were unable to show this. Thus we
were forced to ``layer'' the norm in a certain sense (see section~2
below).

In Section~3 we consider in a wider context $SP_\omega(X)$, the
partially ordered set of all spreading models $(\tilde x_i)$ generated
by weakly null sequences in $X$. The partial order is defined by
domination: we write $(\tilde x_i) \ge (\tilde y_i)$ if for some $C <
\infty$ we have $C \|\sum a_i \tilde x_i \| \ge \|\sum a_i \tilde y_i
\|$ for all scalars $(a_i)$. We identify $(\tilde x_i)$ and $ (\tilde
y_i)$ in $SP_\omega(X)$ if $(\tilde x_i) \ge (\tilde y_i) \ge (\tilde
x_i)$.  We prove (in Proposition~\ref{mainlemma}) that if $C \subseteq
SP_\omega (X)$ is countable then there exists $(\tilde x_i) \in
SP_\omega (X)$ which dominates all members of $C$. This enables us to
prove that in certain cases one can produce an uncountable chain $
\{(\tilde x_i^\alpha)_i\}_{\alpha < {\omega_1}}$ with $(\tilde
x_i^\alpha)_i < (\tilde x_i^\beta)_i$ if $\alpha < \beta < \omega_1$.
The example of the previous section and the above yield a solution to
a uniformity question raised by H.~Rosenthal. The question (and a dual
version) are as follows:
Let a separable Banach space $Z$ have the property that for all
spreading models $(\tilde x_i)$ of normalized basic sequences
$$
\lim_n \left\|\sum^n_{i=1} \tilde x_i\right\|\Big/ n = 0\quad
\left(\text{respectively, } \lim_n \left\|\sum^n_{i=1}\tilde
     x_i\right\| =\infty\right).
$$
Does there exist $(\lambda_n)$ with $\lim\limits_n \lambda_n/n = 0$
(respectively $\lim\limits_n \lambda_n = \infty$) such that for all
spreading models $(\tilde x_i)$ of normalized basic sequences in $Z$
$$
\lim_n \left\|\sum^n_{i=1} \tilde x_i\right\|\Big/\lambda_n = 0
\quad \left(\text{respectively, } \lim_n \left\|\sum^n_{i=1} \tilde
     x_i\right\|\Big/ \lambda_n = \infty\right)?
$$
We give negative answers to these questions.  The example that
solves the first question is the space $X$ of section 2.  Moreover
every subspace of $X$ fails to admit such a sequence $(\lambda_n)$.
We do not know of a hereditary solution to the second question.

In section 4 we consider the problem: if $|SP_\omega(X)| =1$, i.e., if
$X$ has a unique spreading model up to equivalence, must this
spreading model be equivalent to the unit vector basis in $c_0$ or
$\ell_p$ for some $1 \le p < \infty$? The question was asked of us by
S.~A.~Argyros.  It is easy to see that the answer is positive if the
spreading models are uniformly isomorphic. We show that the answer is
positive if 1 belongs to the ``Krivine set'' of some spreading model.

\begin{Def} \label{def1.1}
Let $(x_i)$ be a 1-subsymmetric basic sequence. The {\em Krivine set}
of $(x_i)$ is
the set of $p$'s $(1\le p\le\infty)$ with the following property: \ For all
$\varepsilon>0$ and $n\in\mathbb{N}$ there exists $m\in\mathbb{N}$ and
$(\lambda_k)^m_{k=1} \subset \mathbb{R}$,
such that for all $(a_i)_1^n \subseteq \mathbb{R}$,
$$
\frac1{1+\varepsilon} \|(a_i)^n_{i=1}\|_p \le \left\|\sum^n_{i=1}
   a_iy_i\right\| \le (1+\varepsilon) \|(a_i)^n_{i=1}\|_p
$$
where
$$
y_i = \sum^m_{k=1} \lambda_k x_{(i-1)m+k} \quad \text{for}\quad
i=1,\ldots, n
$$
and $\|\cdot\|_p$ denotes the norm of the space $\ell_p$.
\end{Def}

The proof of Krivine's theorem \cite{K} as modified by Lemberg
\cite{L}) (see also \cite{G}, remark II.5.14 and \cite{MS}), shows
that for every 1-subsymmetric basic sequence $(x_i)$ the Krivine set
of $(x_i)$ is non-empty. It is important to note that our definition
of a Krivine $p$ requires not merely that $\ell_p$ be block finitely
representable in $[x_i]$ but each $\ell_p^n$ unit vector basis is
obtainable by means of an identically distributed block basis.

An immediate consequence of the fact that the Krivine set of a
spreading model is non-empty is the following:

\begin{Rmk}\label{rem1.1}
   Assume that $(x_i)$ is a seminormalized basic sequence in a Banach
   space $X$ which has a spreading model $(\tilde x_i)$.  We can assume
   that for some zero sequence $(\varepsilon_i)$ (\ref{eq0}) is
   satisfied.  Then there is a $p\in[1,\infty]$ such that for all $n$
   and all $\varepsilon>0$ there exists a finite sequence
   $(\lambda_i)_{i=1}^m\subset\R$ so that any block of $(y_i)$ of
   $(x_i)$ of the form
  $$
y_i=\sum_{j=1}^m\lambda_j x_{n(i,j)},
  \text{ with }n(1,1)<n(1,2)<\ldots<n(1,m)<n(2,1)<\ldots
  n(2,m)<n(3,1)\ldots
  $$
  has a spreading model $(\tilde y_i)$ which is isometric to the
  sequence $(\sum_{j=1}^m\lambda_j \tilde x_{(i-1)m+j})_{i\in\N}$ and
  has the property that its first $n$ elements are
  $(1+\varepsilon)$-equivalent to the unit vector basis of $\ell_p^n$.
  For $i_0\in \N$ large enough (or passing to an appropriate
  subsequence of $(x_i)$) we also observe that $(y_{k_j})_{j=1}^n$ is
  $(1+2\varepsilon)$-equivalent to the $\ell_p^n$ unit basis whenever
  $i_0<k_1<\ldots k_n$.
\end{Rmk}

In Section~\ref{sec5} we give sufficient conditions on a Banach space
$X$ for the existence of a subspace $Y$ of $X$ and an operator $T:Y
\to X$ which is not a compact perturbation of the inclusion
map. W.T. Gowers \cite{G} proved that there exists a subspace $Y$ of
the Gowers-Maurey space $GM$ (constructed in \cite{GM}) and there
exists an oparator $T:Y \to GM$ which is not a compact perturbation of
the inclusion map. Here we extend the work of W.T. Gowers to a more
general setting.  For example suppose that $X$ admits a spreading
model $(\tilde x_i)$ which is not equivalent to the unit vector basis
in $\ell_1$ but such that $1$ is in the Krivine set of $(\tilde
x_i)$. Then (Theorem~\ref{Main3}) there exists a subspace $W $ of $ X$
and a bounded operator $T: W \to W$ such that $p(T)$ is not a compact
perturbation of the identity, for any polynomial $p$.

\section{Spreading models containing $\ell_1$ which are not $\ell_1$}
\label{sec2}

\begin{Thm}\label{thm2.1}
   There exists a reflexive Banach space $X$ with an unconditional
   basis such that the spreading model of any normalized basic sequence
   in $X$ is not isomorphic to $c_0$ or $\ell_1$ and is not reflexive.
\end{Thm}

For $x=(x_i)_i \in c_{00}$ we write $\mbox{supp }x= \{ i : x_i \not =
0 \}$. For $x,y \in c_{00}$ and an integer $k$ we say that $x<y$ if
$\max \, \mbox{supp}\, x < \min \, \mbox{supp} \, y$, and we write $k
< x$ if $k < \min \, \mbox{supp}\, x$. $(e_i)$ denotes the unit vector
basis of $c_{00}$.

In order to prove Theorem \ref{thm2.1} we will construct a space $X$
which has certain properties as stated in the following result, which
will easily imply Theorem \ref{thm2.1}.

\begin{Thm}\label{thm2.2}
   There is a space $X$ with the following properties:
\begin{itemize}
\item[a)] $(e_i)$ is a normalized 1-unconditional basis for $X$.

\item[b)] For any normalized block basis of $(e_i)$ having spreading
model $(\tilde x_i)$ we have that $(\tilde x_i)$
is not  equivalent to the unit vector
basis of $\ell_1$.
\item[c)] For any normalized block basic sequence in $X$ having a
   spreading model $(\tilde x_i)$ we have that $\ell_1$ embeds into
   $[(\tilde x_i)]$.
\end{itemize}
\end{Thm}

\begin{proof}[Proof of Theorem \ref{thm2.1}]
   Let $X$ be chosen as in Theorem \ref{thm2.2}.  Since $X$ has an
   unconditional basis and does not contain a subspace isomorphic to
   $\ell_1$ or $c_0$ (otherwise a block basis of $(e_i)$ would be
   equivalent to either the unit vector basis of $\ell_1$ or c$_0$,
   both of which are excluded by (b) and (c)), $X$ must be reflexive.

   Since $X$ is reflexive every normalized basic sequence in $X$ has a
   subsequence which is equivalent to a block basis of $(e_i)$.
   Therefore (b) and (c), and the fact that $\ell_1$ has a unique
   subsymmetric basis, imply that all the spreading models of
   normalized basic sequences in $X$ are neither reflexive nor
   isomorphic to $c_0$ or $\ell_1$.
\end{proof}

\noindent {\bf Construction of the space} $X$: First we choose an
increasing sequence of integers $(n_i)$ such that
\begin{equation}\label{eq1}
\frac1{(n_1+n_2 +\ldots+ n_k)^{1/p}}\sum\limits^k_{i=1} \frac{n_i}{3^i}
\mathop{\longrightarrow}\limits_{k\to\infty} \infty
\quad \text{for all}\quad p>1.
\end{equation}
In order to choose a sequence $(n_i)$ satisfying (\ref{eq1}), first
choose a sequence $(p_k)_k$ with $p_k \searrow 1$ and then inductively
on $k \in \mathbb{N}$ pick $(n_k)$ to satisfy
$$
  \frac1{(n_1+n_2 +\ldots+ n_k)^{1/p_k}} \sum\limits^k_{i=1} 
\frac{n_i}{3^i}    > k
$$
for all $k\in\mathbb{N}$. Now we choose a norm $\|\cdot\|$ on $c_{00}$
to satisfy the
following Tsirelson type equation (see \cite{OS2}):
$$
\|x\| =\|x\|_\infty \vee \sup_{\stackrel{\scriptstyle
     k\in\mathbb{N}, i \le k}{\scriptstyle k\le E^i_1 < E^i_2 <\ldots
     <E^i_{n_i}}} \sum^k_{i=1} \frac1{3^i} \sum^{n_i}_{j=1} \|E^i_jx\|.
$$
Note that we do not require that $E^s_j\cap E^t_{j'}=\emptyset$ if
$s\ne t$. Henceforth in this section $X$ will denote the completion of
$c_{00}$ under this norm. It is easy to see that the unit vector basis
$(e_i)$ is a normalized 1-unconditional basis for $X$. It will be
useful to introduce the sequence of equivalent norms $\|\cdot\|_i$,
for $i\in\mathbb{N}$, as follows:
$$
\|x\|_i = \sup_{E_1 < E_2 <\ldots< E_{n_i}} \sum^{n_i}_{j=1}
\|E_jx\|.
$$
Note that we have
$$
\|x\|  = \|x\|_\infty \vee \sup_{k\in\mathbb{N}} \sum^k_{i=1} \frac1{3^i}
\|[k,\infty)x\|_i.
$$

\begin{proof}[Proof of Theorem \ref{thm2.2}]
a) is immediate.

\noindent
b) We need the following auxiliary results. We
postpone the proofs.

\begin{Lem}\label{lem2.3}
   For any normalized block basis $(y_i)$ of $(e_i)$ and for any
   $\varepsilon>0$ there exists a subsequence $(x_i)$ and
   $i_0\in\mathbb{N}$ such that for any $N\in\mathbb{N}$ and integers
   $k$, $j_1, \ldots , j_N$ with $i_0\le k \le j_1 < j_2< \ldots < j_N$
   we have that
\begin{equation}\label{eq2}
\sum^k_{i=i_0} \frac1{3^i} \left\|[k,\infty) \left(
\frac{1}{N} \sum^N_{s=1} x_{j_s}\right) \right\|_i <
\varepsilon .
\end{equation}
\end{Lem}

\begin{Lem}\label{lem2.4}
   Let $(y_i)$ be a normalized block basis of $(e_i)$ in $X$ which has
   a spreading model $(\tilde y_i)$ and suppose that $N\in\mathbb{N}$
   satisfies
\begin{equation}\label{eq3}
.99 \le \| \frac{1}{2N} ( \tilde y_1 +\ldots+ \tilde y_{2N} )\| .
\end{equation}
Then there exists $k\in\mathbb{N}$ and a subsequence $(x_i)$ of
$(y_i)$ such that for all $j_1 < j_2 <\ldots< j_N$,
\begin{equation}\label{eq4}
.96 < \sum^k_{i=1} \frac1{3^i} \left\|[k,\infty)\left( \frac{1}{N}
\sum^N_{s=1} x_{j_s}\right) \right\|_i .
\end{equation}
\end{Lem}

For the proof of b) assume to the contrary that there exists a
normalized block basis $(y_i)$ of $(e_i)$ whose spreading model
$(\tilde y_i)$ is equivalent to the unit vector basis of $\ell_1$.
Without loss of generality \cite{BL}(Proposition 4 in Chapter II
  Section 2),  we can assume that (\ref{eq3}) is
valid for all $N\in\mathbb{N}$.  For $\varepsilon=.01$ choose
$i_0\in\mathbb{N}$ and a subsequence of $(y_i)$ which satisfies the
conclusion of Lemma~\ref{lem2.3}. Choose $N\in\mathbb{N}$ with
$$
\frac{2}{N} \sum^{i_0-1}_{i=1} n_i < .01 .
$$
Since (\ref{eq3}) is valid, by Lemma~\ref{lem2.4} there exists
$k\in\mathbb{N}$ and a further subsequence $(x_i)$ which satisfies
(\ref{eq4}). Now let $j_1 < j_2 <\ldots< j_N$ with $k\le j_1$ and let
$x = (1/N) \sum^N_{s=1} x_{j_s}$.
We will first estimate for $i\in\N$ the value of $\|x\|_i$. Choose
  $E_1<E_2<\ldots<E_{n_i}$ so that
$$\|x\|_i=\sum_{i=1}^{n_i} \|E_j(x)\|.$$
For $=1,2,\ldots n_i$ put $I_j=\{s\le N: \supp(x_s)\subset E_j\}$ and
$I_0=\{1,2,\ldots, N\}\setminus \bigcup_{j=1}^{n_i} I_j$ and note that
  $I_0=\big\{s\le N:\exists  j_1,j_2\le n_i,j_1\not= j_2,\quad
  \supp(x_s)\cap E_{j_t}\not= \emptyset,\,t=1,2\big\}$,
that $\sum_{j\le n_i} |I_j|\le N$, and that $|I_0|\le \min(N,2n_i)$.
  Moreover note that each $E_j$ can only have a non empty intersection
  with the support of at most two $x_s$'s, $s\in I_0$.
Therefore we deduce
\begin{equation}\label{eq3a}
\|x\|_i=\sum_{j=1}^{n_i} \|E_j(x)\|
       \le\frac1N\sum_{j=1}^{n_i}
  \left[ \sum_{s\in I_j}\|x_s\|+ \Big\|E_j\Big(\sum_{s\in I_0} 
x_s\Big)\Big\|\right]
  \le 1+\min(1,\frac{2n_i}{N}).
\end{equation}

  By Lemma~\ref{lem2.4} we have (the
second term on the right disappears if $k<i_0$)
\begin{align*}
.96  &< \sum^{i_0-1}_{i=1} \frac1{3^i} \|[k,\infty)x\|_i + \sum^k_{i=i_0}
\frac1{3^i} \|[k,\infty)x\|_i\\
&\le \sum^{i_0-1}_{i=1} \frac1{3^i} \|x\|_i + .01  \quad \text{(by
(\ref{eq2}) since } k\le j_1)\\
&\le \sum^{i_0-1}_{i=1} \frac{1}{3^i} \frac{2n_i +N}{N} + .01  \quad
\text{(by (\ref{eq3a}))}\\
&< .01  + .5  + .01= .52
\end{align*}
which is a contradiction.

\noindent c) Here we need the following result whose proof is
again postponed:

\begin{Lem}\label{lem2.5}
   Let $(z_i)$ be a normalized block basis of $(e_i)$ with spreading
   model $(\tilde z_i)$. Then for every $K_1\in\mathbb{N}$ there exists
   $K_2 > K_1$ and $(w_i)$, an identically distributed block basis of
   $(z_i)$, which has a spreading model $(\widetilde w_i)$ (which is
     a block basis of $(\tilde z_i)$) such that for all
   $\ell\in\mathbb{N}\colon \ .98 \le \|w_\ell\|\le 1$ and
\begin{equation}\label{eq5}
\sum^{K_2}_{i=K_1+1} \frac1{3^i} \|[K_2,\infty) w_\ell\|_i >.4.
\end{equation}
\end{Lem}

Let $(z_i)$ be a normalized block basis of $(e_i)$ having a spreading
model $(\tilde z_i)$.
  By passing to a subsequence if necessary we can assume that (\ref{eq0})
  is satisfied for some sequence $(\varepsilon_n)$ which converges to 0.
  By applying Lemma~\ref{lem2.5} repeatedly there
exists an increasing sequence of integers $(K_n)$, $(K_1=0)$, and for
every $n\in\mathbb{N}$ there exists an identically distributed block
basis $(w^{(n)}_i)_i$ of $(z_i)$ having spreading model $(\widetilde
w^{(n)}_i)_i$, which is also a block basis of $(\tilde z_i)$, such
that for all $n,\ell\in \mathbb{N}$, $.98\le \|w^{(n)}_\ell\|\le 1$
and
\begin{equation}\label{eq6}
\sum^{K_{n+1}}_{i=K_n+1} \frac1{3^i}
   \|[K_{n+1},\infty)w^{(n)}_\ell\|_i > .4
\end{equation}
Choose a sequence $(m_i)$ of integers such that $(\widetilde
w^{(i)}_{m_i})$ is a block sequence of $(\tilde z_i)_i$. We claim that
$(\widetilde w^{(i)}_{m_i})$ is equivalent to the unit vector basis of
$\ell_1$. We show that for $(a_i)_{i=1}^N \subseteq \mathbb{R}$,
$$\left\|\sum^N_{i=1} a_i \widetilde w^{(i)}_{m_i}\right\| > .4
\sum_{i=1}^N |a_i|.$$
Let $j_1 < j_2 <\ldots$ be such that
$w^{(1)}_{j_1} < w^{(2)}_{j_2} <\ldots< w^{(N)}_{j_N} <
w^{(1)}_{j_{N+1}} < w^{(2)}_{j_{N+2}} <\ldots< w^{(N)}_{j_{2N}} <
w^{(1)}_{j_{2N+1}} <\ldots$.  Then, since $(z_i)$ satisfies 
(\ref{eq0}), it follows that
$$\left\|\sum^N_{i=1} a_i \widetilde w^{(i)}_{m_i}\right\| = \lim_\ell
\left\|\sum^N_{n=1} a_n w^{(n)}_{j_{(\ell-1)N+n}}\right\|.$$
If we choose $\ell$ such that $\sum^N_{n=1} w^{(n)}_{j_{(\ell-1)N+n}}$
is supported on $[K_{N+1},\infty)$ then
\begin{align}
\label{eq7}
\left\|\sum^N_{n=1} a_n w^{(n)}_{j_{(\ell-1)N+n}}\right\| &\ge
\sum^N_{n=1} |a_n| \sum^{K_{n+1}}_{i=K_n+1} \frac1{3^i} \|[K_{N+1},\infty)
w^{(n)}_{j_{(\ell -1)N+n}}\|_i \\
&= \sum^N_{n=1} |a_n| \sum^{K_{n+1}}_{i=K_n+1} \frac1{3^i}
\|[K_{n+1},\infty) w^{(n)}_{j_{(\ell -1)N+n}}\|_i \nonumber \\
&> .4 \ \sum_{n=1}^N |a_i|  \quad \text{(by (\ref{eq6}))}. \quad \quad
\quad \quad \quad
\quad \quad \quad \quad \quad \quad \quad \quad \quad \quad \quad
\nonumber
\end{align}
\end{proof}

\begin{proof}[Proof of Lemma \ref{lem2.3}]
   Since for all $i$ and $j$ we have $1\le \|y_j\|_i \le n_i$, by a
   simple compactness and diagonalization argument there exists a
   subsequence $(x_i)$ of $(y_i)$ such that
\begin{equation}\label{eq8}
\Big|\|x_i\|_i - \|x_j\|_i\Big| \le 1 \quad \text{for all}\quad i\le j.
\end{equation}
Now we claim that
\begin{equation}\label{eq9}
\sum^\infty_{i=1} \frac1{3^i} \|x_i\|_i \le \frac{3}{2}.
\end{equation}
Indeed, otherwise there exists $k\in \mathbb{N}$ such that
\begin{equation}\label{eq10}
\sum^k_{i=1} \frac1{3^i} \|x_i\|_i > \frac32.
\end{equation}
Choose $j\ge k$ such that $x_j$ is supported on $[k,\infty)$. Then
\begin{align*}
   \|x_j\| &\ge \sum^k_{i=1} \frac1{3^i} \|[k,\infty) x_j\|_i\\
   &= \sum^k_{i=1} \frac1{3^i} \|x_j\|_i \quad \text{(since $x_j$ is
     supported on
     $[k,\infty)$)}\\
   &\ge \sum^k_{i=1} \frac1{3^i} (\|x_i\|_i-1) \quad \text{(by
     (\ref{eq8}) since $j\ge k$)}\\
   &> \frac32 - \sum^k_{i=1} \frac1{3^i} >1 \quad \text{(by
     (\ref{eq10}))}
\end{align*}
which is a contradiction. Thus (\ref{eq9}) is established. Now choose
$i_0\in \mathbb{N}$ such that
\begin{equation}\label{eq11}
\sum^\infty_{i=i_0} \frac1{3^i} \|x_i\|_i
+ \sum^\infty_{i=i_0} \frac1{3^i} < \varepsilon.
  \end{equation}
Let $k, j_1, \ldots , j_N \in\mathbb{N}$ with $i_0\le k \le j_1 < j_2
<\ldots<  j_N$. We have
\begin{align*}
   \sum^k_{i=i_0} \frac1{3^i} \left\|[k,\infty) \sum^N_{s=1}
     x_{j_s}\right\|_i
   &\le \sum^N_{s=1} \sum^k_{i=i_0} \frac1{3^i} \|x_{j_s}\|_i\\
   &\le \sum^N_{s=1} \sum^k_{i=i_0} \frac1{3^i} (\|x_i\|_i+ 1) \quad
   \text{(by
     (\ref{eq8}) since $k\le j_1$)}\\
   &< N\varepsilon \quad \text{(by (\ref{eq11}))}.
\end{align*}
\end{proof}

\begin{proof}[Proof of Lemma \ref{lem2.4}]
   From (\ref{eq3}) there exists a subsequence $(z_i)$ of $(y_i)$ such
   that
\begin{equation}\label{eq12}
.98 < \| \frac{1}{2N} ( z_1+z_2 +\ldots+ z_N + z_{j_1} +z_{j_2}
+\ldots+  z_{j_N}) \|
\end{equation}
for all $N < j_1 < j_2 <\ldots < j_N$. Let $K$ be the maximum element
in the support of $z_N$. Now for $j_1 < j_2 <\ldots< j_N$ let $u =
(z_1 +\ldots+ z_N)/N$, $v = (z_{j_1} +\ldots+ z_{j_N})/N$ and
$w=(u+v)/2$. By the definition of the norm of $X$ there exists $k'\in
\mathbb{N}$, which depends on $j_1,\ldots, j_N$, such that
$$\|w\| = \sum^{k'}_{i=1} \frac1{3^i} \|[k',\infty)w\|_i .$$
By (\ref{eq12}) we have that $.98 < \|w\| $ and thus $k'\le K$.
By the triangle inequality we obtain
\begin{align*}
   .98 &< \frac{1}{2}\sum^{k'}_{i=1} \frac1{3^i} \|[k',\infty)u\|_i +
   \frac12 \sum^{k'}_{i=1}
   \frac1{3^i} \|[k',\infty)v\|_i\\
   &\le \frac12 \|u\| + \frac12 \sum^{k'}_{i=1} \frac1{3^i}
   \|[k',\infty)v\|_i \le .5 +\frac12 \sum^{k'}_{i=1} \frac1{3^i}
   \|[k',\infty)v\|_i.
\end{align*}
Thus
$$.96 < \sum^{k'}_{i=1} \frac1{3^i} \|[k',\infty)v\|_i.$$
Now by Ramsey's theorem \cite{Ra} (see also \cite{O}) there exists a
subsequence $(x_i)$ of $(z_{N+i})$ and $k\le K$ such that (\ref{eq4})
is valid for all $j_1 < j_2 <\ldots <j_N$.
\end{proof}

\begin{proof}[Proof of Lemma \ref{lem2.5}]
   Let us first note that neither $\ell_p$, $p>1$, nor $c_0$ are
   finitely block represented in $X$.  Indeed, if $(x_i)$ for
   $i=1,\ldots, n_1+\ldots+ n_k$ (some $k\in \mathbb{N}$) is a
   normalized block basis of $(e_i)$ which is 2-equivalent to the first
   $n_1+\ldots+ n_k$ unit basic vectors of $\ell_p$ for some $p>1$,
   then if $\mbox{supp }x_1>k$, it follows that
\begin{align*}
   2(n_1+\ldots+n_k)^{1/p} &\ge \left\| \sum^{n_1
       +\ldots+ n_k}_{i=1} x_i\right\|\\
   &\ge \sum^k_{i=1} \frac1{3^i} \sum^{n_i}_{j=n_{i-1}+1} \|x_j\| \quad
   \text{(by definition of the norm; set $n_0=0$)}\\
   &= \sum^k_{i=1} \frac1{3^i} n_i
\end{align*}
which contradicts (\ref{eq1}). Similarly the case $p=\infty$ is
excluded and thus the conclusions of  Remark~\ref{rem1.1} hold only
for $p=1$.

Let $(z_i)$ be a normalized block sequence in $X$ having a spreading
model $(\tilde z_i)$, and let $K_1\in \mathbb{N}$. Choose
$N\in\mathbb{N}$ such that
\begin{equation}\label{eq13}
\frac2N \sum^{K_1}_{i=1} n_i < .01.
\end{equation}
By Remark~\ref{rem1.1} there exists an identically distributed block
basis $(y_i)$ of $(z_i)$ having spreading model $(\tilde y_i)$ which
satisfies (\ref{eq3}) and ($\tilde y_\ell)$ is a block basis of
$(\tilde z_i)$. Thus by Lemma~\ref{lem2.4} there exists
$K_2\in\mathbb{N}$ and a subsequence $(x_i)$ of $(y_i)$ such that
(\ref{eq4}) is satisfied for $k=K_2$ and for all $j_1 < j_2
<\ldots<j_N$. Let
$$w_\ell = \frac1N \sum^N_{j=1} x_{N(\ell -1)+j}\quad \text{for}\quad
\ell\in\mathbb{N}.$$
Since (\ref{eq3}) is satisfied, by passing to a
subsequence we can assume that $.98 \le \|w_\ell\| \le 1$ for all
$\ell$. Let $(\tilde w_i)$ be the spreading model of $(w_i)$. Then for
all $\ell\in\mathbb{N}$,
$$\widetilde w_\ell = \frac1N \sum^N_{j=1} \tilde y_{N(\ell -1)+j}.$$
Thus $(\widetilde w_\ell)$ is a block basis of $(\tilde z_i)$ and
\begin{equation}\label{eq14}
.96 < \sum^{K_2}_{i=1} \frac1{3^i} \|[K_2,\infty) w_\ell\|_i.
\end{equation}
Note also that by  with the same argument as in the proof of (\ref{eq3a}),
\begin{equation}\label{eq15}
\sum^{K_1-1}_{i=1} \frac1{3^i} \|[K_2,\infty)w_\ell\|_i\le
\sum^{K_1-1}_{i=1} \frac1{3^i} \frac{2n_i +N}{N}
  < .01 + .5= .51.\quad \text{(by (\ref{eq13}))}
\end{equation}
Now (\ref{eq14}) and (\ref{eq15}) immediately give (\ref{eq6}).
\end{proof}

\section{The set of spreading models of $X$}
\label{sec3}

We recall the standard

\begin{Def}\label{D3.1}
   Let $(x_i)$ and $(y_i)$ be basic sequences and $C \ge 1$. We say
   that $(x_i)$ $C$-dominates $(y_i)$, if $C \| \sum_i a_i x_i \| \ge
   \| \sum a_i y_i \|$ for all $(a_i) \in c_{00}$. We say that $(x_i)$
   dominates $(y_i)$, denoted by $(x_i) \ge (y_i)$, if $(x_i)$
   $C$-dominates $(y_i)$ for some $C \ge 1$.
   We write $(x_i) > (y_i)$,
   if $(x_i) \ge (y_i)$ and $(y_i) \not\ge (x_i)$.
   If $(x^n_i)_i$ (for $n \in
   \N$) is a sequence of basic sequences and $(z_i)$ is a basic
   sequence, then we say that $(z_i)$ uniformly dominates the sequences
   $(x^n_i)_i$ if there exists $C \ge 1$ such that $(z_i)$
   $C$-dominates $(x^n_i)_i$ for all $n \in \N$.
\end{Def}

The set $SP(X)$ of all spreading models generated by normalized basic
sequences in $X$ is partially ordered by domination, provided that we
identify equivalent spreading models. $SP_\omega(X)$ denotes the
subset of those spreading models generated by weakly null sequences.

Our first result in this section shows that every
countable subset of $SP_\omega(X)$ admits an upper bound in
$SP_\omega(X)$.

\begin{Prop} \label{mainlemma}
Let $(C_n) \subset (0, \infty)$ such that $\sum C_n^{-1} < \infty$ and
for $n \in \N$ let $(x^n_i)_i$ be a normalized weakly null sequence
in some Banach space $X$ having spreading model
$(\tilde{x}^n_i)_i$. Then there exists a seminormalized weakly null
basic sequence $(y_i)$ in $X$
  with spreading model $(\tilde y_i)$ such that $(\tilde y_i)$ 
$C_n$-dominates $(\tilde x^n_i)_i$
for all $n \in \N$. \\
Furthermore, if
for all $n \in \N$, $(\tilde{x}^n_i)_i$ is not equivalent to the unit
vector basis of $\ell_1$, then $(\tilde{y}_i)$ is not equivalent to 
the unit vector
basis of $\ell_1$.\\
Moreover, if $(z_i)$
is a basic sequence which uniformly dominates $(\tilde{x}^n_i)_i$ for
all $n \in \N$, then $(z_i)$ dominates $(\tilde y_i)$. \\
\end{Prop}

In order to prove Proposition~\ref{mainlemma} we need first to
generalize the fact that spreading models
of normalized weakly null sequences exist and are suppression
1-unconditional.

Lemma \ref{lem3.3b} is actually a special case of a
more general situation \cite{HO}. The results could also
be phrased in terms of countably branching trees of order $mn$
  and proved much like the arguments in \cite{KOS}.

\begin{Lem} \label{lem3.3b}
   Let $n, m \in \mathbb{N}$ and $\varepsilon >0$. Let $(x^{(1)}_i)_i,
   (x^{(2)}_i)_i, \ldots , (x^{(n)}_i)_i$ be normalized weakly null
   sequences in a Banach space $X$.  Then there exists a subsequence
   $L$ of $\N$ so that for all families of integers
   $(k^{(i)}_j)_{i=1,j=1}^{n,m}$ and $(\ell^{(i)}_j)_{i=1,j=1}^{n,m}$ 
in $L$, with $k_1^{(1)}<k_1^{(2)}<
   \ldots k_1^{(n)}<k_2^{(1)}< \ldots < k_2^{(n)}< \ldots <k_m^{(1)}<
   \ldots <k_m^{(n)}$ and $\ell_1^{(1)} <\ell_1^{(2)} < \ldots <
   \ell_1^{(n)} < \ell_2^{(1)}< \ldots < \ell_m^{(n)}$, and
   $(a^{(j)}_i)_{i=1,j=1}^{m,n} \subseteq [-1,1]$ we have
$$
\left|\bigl\| \sum_{i=1}^m \sum_{j=1}^n a^{(j)}_i x^{(j)}_
{\ell_i^{(j)}} \| -
\| \sum_{i=1}^m \sum_{j=1}^n a^{(j)}_i x^{(j)}_{k_i^{(j)}}
\bigr\|\right|
\le \varepsilon. $$
\end{Lem}

\begin{proof} This follows easily by Ramsey's theorem.
   Let $(a^{(j)}_i)_{i=1,j=1}^{m,n} \subseteq [-1,1]$ such that not all
   the $a^{(j)}_i$'s are $0$.  Partition $[0,mn]$ into finitely many
   intervals of length less $\varepsilon/2$.  Partition the finite
   sequences $k_1^{(1)}<k_1^{(2)}< \ldots < k_m^{(n)}$ of $\mathbb{N}$
   according to which interval $\| \sum_{i=1}^m \sum_{j=1}^n a^{(j)}_i
   x^{(j)}_{k_i^{(j)}} \|$ belongs.  Thus by Ramsey's theorem for some
   infinite subsequence $L$ of $\N$ these expressions belong to the
   same interval, if $k_i^{(j)} \in L$ for $i=1, \ldots , m$ and $j=1,
   \ldots , n$. We repeat this for a finite $\varepsilon/4$-net of
   $[-1,1]^{mn}$ endowed with the $\ell_1^{mn}$ norm. \end{proof}

\begin{Lem} \label{lem3.3d}
   Let $n,m \in \mathbb{N}$, $\varepsilon >0$ and $(x^{(1)}_i)_i,
   \ldots , (x^{(n)}_i)_i$ be normalized weakly null sequences in a
   Banach space $X$. Then there exists a subsequence $L$ of $\N$ so
   that for all integers in $L$, $k_1^{(1)}<k_1^{(2)}< \ldots
   k_1^{(n)}<k_2^{(1)}< \ldots < k_2^{(n)}< \ldots <k_m^{(1)}< \ldots
   <k_m^{(n)}$, the vectors $(y^{(j)}_{k^{(j)}_i})_{i=1,j=1}^{m,n}$
   form a suppression $(1+\varepsilon)$-unconditional basic sequence.
\end{Lem}

\begin{proof}
   By passing to subsequences, if necessary, we may assume that the sequence
   $(x^{(j)}_i)_{j=1,i=1}^{n,\infty}$ satisfies the conclusion of
   Lemma~\ref{lem3.3b} for $\varepsilon$ replaced by $\varepsilon /2$ 
and $L=\N$.
   Let $\delta = \varepsilon /(2nm)$. For every $j_0 \leq n$, since
   $(x^{(j_0)}_i)_i$ is weakly null we have that for every $i_0 \in
   \mathbb{N}$ there exists $i_1 >i_0$ such that for every functional
   $f \in X^*$ of norm $1$ there exists $i \in [i_0,i_1]$ with
   $|f(x^{(j_0)}_i)|< \delta$. Iterating this procedure we can pass to
   an infinite subsequence $L$ of $\mathbb{N}$ with the following
   property: Given $k_1^{(1)}<k_1^{(2)}< \ldots < k_1^{(n)}< k_2^{(1)}<
   \ldots k_m^{(n)}$ in $ L$, $F \subseteq \{ k_1^{(1)}, k_1^{(2)},
   \ldots , k_1^{(n)}, k_2^{(1)}, \ldots , k_m^{(n)} \}$, and $f \in
   X^*$ of norm $1$, there exist $\ell_1^{(1)}< \ell_1^{(2)}< \ldots <
   \ell_1^{(n)}< \ell_2^{(1)}< \ldots < \ell_m^{(n)}$ in $\mathbb{N}$
   with $\ell_i^{(j)} = k_i^{(j)}$ if $k_i^{(j)} \not \in F$ and
   $|f(x^{(j)}_{\ell_i^{(j)}})| < \delta$ if $k_i^{(j)} \in F$. Let
   $k_1^{(1)}<k_1^{(2)}< \ldots < k_1^{(n)}< k_2^{(1)}< \ldots <
   k_m^{(n)}$ in $L$, $F \subseteq \{ k_1^{(1)}, k_1^{(2)}, \ldots ,
   k_1^{(n)} , k_2^{(1)}, \ldots , k_m^{(n)} \}$, and
   $(a^{(j)}_i)_{i=1,j=1}^{m,n} \subseteq [-1,1]$ with
   $\|\sum_{i=1}^m\sum_{j=1}^n a^{(j)}_i x^{(j)}_{k^{(j)}_i}\|=1$.
   There exists $f \in X^*$ of norm $1$ such that
\begin{align*}
   \| \sum_{ \{ (i,j): k_i^{(j)} \not \in F \} } a^{(j)}_i
   x^{(j)}_{k_i^{(j)}} \| & = f\Bigl(\sum_
   { \{ (i,j): k_i \not \in F \} }
   a^{(j)}_i x^{(j)}_{k_i^{(j)}} \Bigr), \text{ and choosing
$(\ell^{(j)}_i)$ as above, }\\
   & \leq f\Bigl( \sum_{i=1}^m \sum_{j=1}^n a^{(j)}_i
   x^{(j)}_{\ell^{(j)}_i} \Bigr) +  \delta n m  \\
   & \leq \| \sum_{i=1}^m \sum_{j=1}^n a^{(j)}_i x^{(j)}_{\ell_i^{(j)}}
   \| + \frac{\varepsilon}{2} \leq (1 + \frac{\varepsilon}{2} )
   +\frac{\varepsilon}{2} = 1+\varepsilon .
\end{align*}
\end{proof}

\begin{proof}[Proof of Proposition~\ref{mainlemma}]
   Using Lemma \ref{lem3.3d}, a diagonal argument and relabeling we can
   assume that for all $k$ the sequences $(x^1_i)_i, \ldots ,
   (x^k_i)_i$ satisfy the conclusion of Lemma \ref{lem3.3d} for $m\le
   n=k$, $\varepsilon =1$ provided that $k \le k^{(1)}_1$. Let $m_1=0$
   and for $i \in \N$ let $m_{i+1}= m_i + i$. Let $(C_i) \subset (0,
   \infty)$ such that $\sum C_i^{-1} < \infty$. By passing to
   subsequences of $(x^{n}_i)_i$, for each $n\in\N$,
  we can assume in addition that $(y_j)$,
  where
$$
y_j = \sum_{i=1}^j 16 C_i^{-1} x^i_{m_j +i} \mbox{ for all }j\in\N,
$$
is is a seminormalized block  and has a spreading model
  $(\tilde y_j)$.
It is easy to check that $(y_j)$ is weakly null since each
$(x^n_i)_{i=1}^{\infty}$ is  weakly null.
Let $i_0, m \in \N$ and  $(a_j)_{j=1}^m \subset \R$. Let $n
\in \N$ such that $m \le n$ and
$$
8 \sum_{j=1}^m | a_j | \sum_{i=n+1}^\infty C_i^{-1} \le
C_{i_0}^{-1} \| \sum_{j=1}^m a_j \tilde{x}^{i_0}_j \| .
$$
Let $k_1<k_2< \cdots <k_m$ be such that $n <k_1$, and in addition, for
$k_1 \le \ell_1< \cdots < \ell_m$,
$\frac{1}{2} \| \sum_{j=1}^m a_j x^{i_0} _{\ell_j} \|
\le \| \sum_{j=1}^m a_j \tilde{x}^{i_0}_j \|
\le 2 \| \sum_{j=1}^m a_j x^{i_0} _{\ell_j} \|$ for all $i=1, \ldots , n$,
and $\frac{1}{2} \| \sum_{j=1}^m a_j y_{k_j}\| \le \| \sum_{j=1}^m a_j
\tilde{y}_j \| \le 2\| \sum_{j=1}^m a_j y_{k_j}\| $.
We have from Lemma \ref{lem3.3d} and our inequalities,
\begin{align}\label{B*}
\| \sum_{j=1}^m a_j \tilde{y}_j \| \ge & \frac{1}{2} \| \sum_{j=1}^m
a_j y_{k_j}\| \\
\ge & \frac{1}{2} \| \sum_{j=1}^m \sum_{i=1}^n a_j 16 C_i^{-1} 
x^i_{m_{k_j} +i} \|
- \frac{1}{2} \sum_{j=1}^m \sum_{i=n+1}^\infty |a_j| 16 C_i^{-1} \notag  \\
\ge & \frac{1}{4} \| \sum_{j=1}^m a_j 16 C_{i_0}^{-1}x^{i_0}_{m_{k_j} +{i_0}}\|
- C_{i_0}^{-1} \| \sum_{j=1}^m a_j \tilde{x}^{i_0}_j \|\notag \\
\ge & \frac{1}{8} \| \sum_{j=1}^m a_j 16 C_{i_0}^{-1}\tilde{x}^{i_0}_j \|
- C_{i_0}^{-1} \| \sum_{j=1}^m a_j \tilde{x}^{i_0}_j \| \notag\\
= & C_{i_0}^{-1} \| \sum_{j=1}^m a_j \tilde{x}^{i_0}_j \| .\notag
\end{align}
This proves the first part of the proposition. In order to show the remaining
  part let $m\in\N$ and $(a_j)_{j=1}^m
\subseteq \R$, and note that
\begin{align} \label{A*}
\| \sum_{j=1}^m a_j \tilde{y}_j \|
&= \lim_{j_1\to\infty}\ldots \lim_{j_m\to\infty}
\Big\|\sum_{s=1}^m a_s y_{j_s}\Big\|\\
&= \lim_{j_1\to\infty}\ldots \lim_{j_m\to\infty}
\Big\|\sum_{s=1}^m a_s \sum_{i=1}^{j_s} 16 C_i^{-1} 
x^i_{m_{j_s}+i}\Big\| \notag\\
&\le \lim_{j_1\to\infty}\ldots \lim_{j_m\to\infty}
   \sum_{i=1}^\infty 16 C^{-1}_i\Big\|\sum_{s=1}^m a_s 
x^i_{m_{j_s}+i}\Big\|\notag\\
&=\sum_{i=1}^\infty 16 C^{-1}_i\Big\|\sum_{s=1}^m a_s  \tilde 
x^i_s\Big\|. \notag
\end{align}
Now the second part of the claim follows from (\ref{B*}) and (\ref{A*})
\end{proof}

\begin{Rmk}\label{increas}Using a similar argument we can prove the following:
   \begin{itemize}
   \item [1)]
  Let $C =\{(\tilde x_i^n)_i\}_{n \in \N}$ be a strictly increasing chain in
$SP_\omega(X)$. Suppose that $(\tilde z_i)_i \in SP_\omega (X)$ is an
upper bound for $C$. Then there exists an upper bound $(\tilde x_i)_i
\in SP_\omega (X)$ for $C$ which is not equivalent to $(\tilde
z_i)_i$.

\item [2)] %% The argument also yields that
If $(\tilde x_i^n)_i \in
   SP_\omega (X)$ for $n \le m \in \N$ then there exists
$(\tilde x_i)_i \in SP_\omega (X)$ which is equivalent to the norm
given by
$$
\| (a_i)\| = \max_{n \le m} \|\sum_i a_i \tilde x_i^n\|.
$$
An analogous result for asymptotic structure of spaces with a
shrinking basis is obtained in \cite{MT} (Proposition 5.1).
   \end{itemize}
\end{Rmk}

\bigskip

\begin{Prop}\label{pro3.3}
  Suppose that $(x_i)$ is a normalized weakly null sequence in a Banach
  space $X$ which has a spreading model $(\tilde x_i)$ which is not
  equivalent to the unit vector basis of $\ell_1$.  Assume that $1$
  belongs to the Krivine set of $(\tilde x_i)$.  Then for all sequences
  $(\lambda_n)\subset\R$, with $\lambda_n\nearrow\infty$ and $\lim_n
  \frac{n}{\lambda_n}=\infty$ there is a block sequence $(y_n)$ of
  $(x_n)$ having a spreading model $(\tilde y_n)$ which satisfies:
\begin{equation*}
\limsup_n\frac{n}{\|\sum_{i=1}^n \tilde y_i\|}=
\limsup_n\frac{\|\sum_{i=1}^n \tilde y_i\|}{\lambda_n}=\infty.
\end{equation*}
Moreover, the set of all spreading models in $X$ generated by weakly
null sequences and not equivalent to $\ell_1$ has no maximal element
(with respect to domination).
\end{Prop}

Note that the space $X$ constructed in Section~2 is reflexive and satisfies
the hypothesis of the proposition (as does every subspace of $X$).
% , hence
% Question~\ref{qst3.1} has a strongly negative  answer.

\begin{proof}%%[Proof of Proposition~\ref{pro3.3}]
Using $\lim
\frac{n}{\lambda_n}= \infty$, choose a subsequence $(n_k)$ of $\N$ such
that $n_k / \lambda_{n_k} \ge 2^{k+1} k$ for all $k$. Since $1$ belongs to
the Krivine set of $(\tilde{x}_i)$,
for every $n \in \N$ there exists a block sequence $(x^n_i)_i$ of $(x_i)$
which is identically distributed with respect to $(x_i)$ and it has
spreading model $(\tilde{x}^n_i)_i$ as given in Remark~\ref{rem1.1}
(for $p=1$ and $\varepsilon=1$) satisfying
$$
\frac{n_k}{2} \le \| \sum_{i=1}^{n_k} \tilde{x}^k_i \| .
$$
Since $(x_i)$ is weakly null and not
equivalent to the unit vector basis of $\ell_1$ we have that for all
$n \in \N$, $(x^n_i)_i$ is weakly null and $(\tilde{x}^n_i)_i$ is not
equivalent to the unit vector basis of $\ell_1$. We can also assume
without loss of generality that $(x^n_i)_i$ is normalized. Let
$(y_i)$ be the sequence which is provided by Proposition~\ref{mainlemma} for
$C_k=2^{-k}$. By
the ``furthermore'' part of Proposition~\ref{mainlemma} we have that
$(\tilde{y}_i)$ is not equivalent to the unit vector basis of $\ell_1$
thus
$$
\limsup_n \frac{n}{\| \sum_{i=1}^n \tilde{y}_i \| } = \infty .
$$
Also, by Proposition~\ref{mainlemma} we have that
$$
2^k \| \sum_{j=1}^{n_k} \tilde{y}_j \| \ge \| \sum_{j=1}^{n_k}
\tilde{x}^k_i \| \ge \frac{n_k}{2} .
$$
Thus for all $k \in \N$,
$$
\frac{\| \sum_{j=1}^{n_k} \tilde{y}_j \| }{\lambda_{n_k} } \ge
\frac{n_k}{2^{k+1} \lambda_{n_k}} \ge k,
$$
which shows that
$$
\limsup_n \frac{ \| \sum_{j=1}^n \tilde{y}_j \| }{ \lambda_n } =
\infty ,
$$
and finishes the proof of the first part of Proposition~\ref{pro3.3}.

To prove the moreover part, given a spreading model
$(\tilde{z}_i) \in SP_\omega(X) $ not equivalent to $\ell_1$
use the  first part of the Proposition to get $(\tilde{y}_i)$
with
$$
\limsup_n\frac{\|\sum_{i=1}^n \tilde y_i\|}{\|\sum_{i=1}^n \tilde z_i\|
}=\infty.
$$
Then any spreading model  $(\tilde w_i)$ which dominates
both $(\tilde{z}_i)$ and $(\tilde{y}_i)$ (and  exists
by Proposition~\ref{mainlemma}) is not equivalent to
$(\tilde{z}_i)$, which therefore  is  not maximal.
\end{proof}

In some circumstances we will be able to conclude that $SP_\omega(X)$
admits a transfinite strictly increasing chain. The logical part of
the argument is a simple proposition.

\begin{Prop}\label{general}
Let $X$ be a separable infinite dimensional  Banach space.
Let $C \subseteq SP_\omega(X)$ be a non-empty set satisfying
the following two conditions:
\begin{itemize}
\item[(i)]
$C$ does not have a maximal element with respect to domination;
\item[(ii)]
for every
$(\tilde{X}_n) \subseteq C$ there exists $\tilde{X} \in C$
such that $\tilde{X}_n \le \tilde{X}$ for every $n \in \N$.
\end{itemize}
Then for all $\alpha < \omega_1$ there exists $\tilde{X}^\alpha
\in C$ such that if $\alpha < \beta < \omega_1$ then
$\tilde{X}^\alpha < \tilde{X}^\beta$.
\end{Prop}

\begin{proof}
We use transfinite induction. Suppose that
$\tilde{X}^\alpha $ have been constructed for
$\alpha < \beta < \omega_1$. Then $\tilde{X}^\beta $
is chosen using (i) if $\beta$ is a successor ordinal and (ii)
if $\beta$ is a limit ordinal.
\end{proof}

\begin{Rmk}\label{nnew}
(1)
  The set $C = SP_\omega(X)$ satisfies condition (ii) by
virtue of Proposition~\ref{mainlemma}. Hence if $SP_\omega(X)$ does
not have a maximal element, then it contains an uncountable
increasing chain.

(2)
Suppose $SP_\omega(X)$ contains $(\tilde x_i)$ such that $1$ is in
the Krivine set of $(\tilde{x}_i)$ but $(\tilde{x}_i)$ is not
equivalent to the unit vector basis of $\ell_1$. Let $C$ be the set of
all elements of $SP_\omega(X)$ which are not equivalent to the unit
vector basis in $\ell_1$. Then it satisfies (ii) by
Proposition~\ref{mainlemma} and (i) by
Proposition~\ref{pro3.3}. Therefore $C$ contains an uncountable
increasing chain. Examples of such a space $X$ are the space
constructed in Section 2,
Gowers-Maurey space $GM$ (\cite{GM}) and
Schlumprecht's space $S$  (\cite{S2}).
\end{Rmk}

\bigskip

\bigskip

The following result is a strengthening of Proposition
\ref{pro3.3}.
First recall that if
$(x_i)_{i=1}^n$ and $(y_i)_{i=1}^n$ ($n \in \mathbb{N}$) are two basic
sequences then the basis-distance between them is
defined by
$$
d_b \left( (x_i)_{i=1}^n, (y_i)_{i=1}^n \right) = \sup \left\{
\frac{ \| \sum_{i=1}^n a_i x_i \| }{\| \sum_{i=1}^n b_i x_i \|}:
\| \sum_{i=1}^n a_i y_i \| = \| \sum_{i=1}^n b_i y_i \| =1 \right\} .
$$

\begin{Prop} \label{pro3.3a}
   Let $(z_i)$ be a normalized basis and $C<\infty$. Let $X$ be an
   infinite dimensional Banach space such that the spreading model
   $(\tilde x_i)$ of any normalized weakly null sequence $(x_i)$ in $X$
   is $C$- dominated by $(z_i)$.  Assume also that for all $n \in
   \mathbb{N}$ there exists a normalized weakly null sequence
   $(x_i^{n})_i$ in $X$ with spreading model $(\tilde{x}^{n}_i)_i$ such
   that $(\tilde{x}_i^n)_{i=1}^n$ $C$-dominates $(z_i)_{i=1}^n$ for all
   $n \in \N$.  Then for every $\lambda_n \nearrow \infty$ there exists
   a normalized weakly null sequence $(y_i)$ in $X$ with spreading
   model $(\tilde{y}_i)$ so that
$$
\liminf_n \frac{d_b \left( (\tilde{y}_i)_{i=1}^n, (z_i)_{i=1}^n
\right)}{\lambda_n} =0.
$$
\end{Prop}

\begin{proof}
Since $\lambda_n \nearrow \infty$, choose a sequence $(n_k)$ of
integers such that $k 2^k \le \lambda_{n_k}$ for all $k$. Apply
Proposition~\ref{mainlemma} to obtain a seminormalized weakly null sequence
$(y_i)$ in $X$ with a spreading model $(\tilde{y}_i)$ such that
$(\tilde{y}_i)$ $2^k$-dominates $(\tilde{x}^{n_k}_i)$, for all $k \in
\N$. By the moreover part of Proposition~\ref{mainlemma} we also have that
there exists $C' < \infty$ such that $(z_i)$ $C'$-dominates
$(\tilde{y}_i)$. Let $k \in \N$ and $(b_i)_{i=1}^{n_k}$ be a sequence
of scalars. Then
$$
\| \sum_{i=1}^{n_k} b_i \tilde{y}_i \| \ge 2^{-k} \| \sum_{i=1}^{n_k}
b_i \tilde{x}^{n_k}_i \| \ge 2^{-k}C^{-1} \| \sum_{i=1}^{n_k} b_i z_i
\| .
$$
Thus for $k \in \N$, if $(a_i)_{i=1}^{n_k}$ and $(b_i)_{i=1}^{n_k}$
are finite sequences of scalars satisfying
$\| \sum_{i=1}^{n_k} a_i z_i \| = \| \sum_{i=1}^{n_k} b_i z_i \| =1$, then
$$
\frac{\| \sum_{i=1}^{n_k} a_i \tilde{y}_i \| }{\| \sum_{i=1}^{n_k} 
b_i \tilde{y}_i \| }
\le \frac{C'}{2^{-k}C^{-1}} = C C' 2^k \le \frac{ \lambda_{n_k}}{k}CC'.
$$
Hence
$d_b \left( (\tilde{y}_i)_{i=1}^{n_k}, (z_i)_{i=1}^{n_k} \right)
/\lambda_{n_k} \le k^{-1} CC' $
which tends to zero. The result follows by normalizing $(y_i)$.
\end{proof}

Propositions \ref{pro3.3} and \ref{pro3.3a} motivate the following

\begin{Qst} \label{qst3.11}
Which normalized subsymmetric bases $(y_i)$ (if any) have the
following property: If $X$ is a separable infinite dimensional Banach
space so that no spreading model of $X$ is equivalent to $(y_i)$ then
there exists $\lambda_n \nearrow \infty$ and a subspace $Y$ of $X$
such that for all spreading
models $(\tilde x_i)$ of normalized basic sequences in $Y$,
$$
\liminf_j \frac{d_b \left( (\tilde{x}_i)_{i=1}^j, (y_i)_{i=1}^j
\right)}{\lambda_j} > 0.
$$
\end{Qst}

This question is a generalization of the following problem raised by
H.~Rosenthal (which is solved by Proposition~\ref{pro3.3})

\begin{Qst}\label{qst3.1}
Let $Z$ be a separable infinite dimensional Banach space so that whenever
$(\tilde x_i)$ is the spreading model of a normalized basic sequence in $Z$
then
$$\lim_n \left\|\sum^n_{i=1} \tilde x_i\right\|\Big/ n = 0.$$
Does there exist
$\lambda_n\nearrow \infty$ such that  $\lim_n \lambda_n/n = 0$ and
for all spreading models $(\tilde x_i)$ of normalized basic
sequences in $Z$
$$\lim_n \left\|\sum^n_{i=1} \tilde x_i\right\|\Big/ \lambda_n = 0?$$
\end{Qst}

The hypothesis on $Z$ is equivalent to: $Z$ is reflexive and no
spreading model of $Z$ is isomorphic to $\ell_1$. The question asks
whether all spreading models of $Z$ must be uniformly distancing
themselves from $\ell_1$ for large enough dimensions.

Question~\ref{qst3.1} just asks if one could take $(y_i)$
in Question~\ref{qst3.11} to be the unit vector basis of  $\ell_1$.
Proposition~\ref{pro3.3} shows that this is not true, even hereditarily.

The version of  Question~\ref{qst3.11} for the unit vector basis of $c_0$ is
the following

\begin{Qst}\label{qst3.2}
Let $Z$ be a separable infinite dimensional Banach space so that whenever
$(\tilde x_i)$ is a spreading model of a normalized basic sequence in $Z$ then
$$\lim_n \left\|\sum^n_{i=1} \tilde x_i\right\| = \infty.$$
Does there exist
a sequence $(\lambda_n)$ with $\lambda_n\nearrow \infty$ such
that for all spreading models $(\tilde x_i)$ of normalized basic
sequence in $Z$
$$\lim_n \left\|\sum^n_{i=1} \tilde x_i\right\|\Big/ \lambda_n = \infty?$$
\end{Qst}

The hypothesis of this question is equivalent to: no spreading model
of $Z$ is isomorphic to $c_0$.

We now construct a space $Y$ which answers  Question~\ref{qst3.2}
in the negative.

  First
recall the notion of the dyadic tree $T$ and the partial order $\prec$ on it:
$$T = \{(n,i)\colon \ n=0,1,2,\ldots, 0\le i < 2^n\}.$$
We define $(n,i) \prec (m,j)$ if and only if $n<m$ and there exist integers
$i = i_0 < i_1 <\ldots <i_k =j$ $(k=m-n)$ such that $i_1 \in \{2i,2i+1\}$,
$i_2\in \{2i_1, 2i_1+1\},\ldots, i_k \in \{2i_{k-1}, 2i_{k-1}+1\}$. We also
write $(n,i)\preceq (m,j)$ if either $(n,i)\prec (m,j)$ or $(n,i) = (m,j)$. A
branch of $T$ will be a subset of $T$ of the form $\{(0,0), (1,i_1),
(2,i_2),\ldots\}$ with $i_k\in \{2i_{k-1}, 2i_{k-1}+1\}$ for all $k$. Let
$\mathcal{B}$ denote the set of all branches. We can naturally identify
$\mathcal{B}$ with the set of infinite sequences of 0 and 1's. Thus we can also
identify $\mathcal{B}$ with the set $\mathcal{B}'$ of finite or infinite
sequences of positive integers (the empty sequence of $\mathcal{B}'$
corresponds to $(0,0,\ldots)$; the non-empty sequence $(b_n)\in \mathcal{B}'$
corresponds to the sequence of 0 and 1's where the 1's are located at the
positions $b_1,b_1+b_2$, $b_1+b_2+b_3,\ldots$). Every branch $b\in\mathcal{B}$
corresponds in this bijective way to a $b'\in\mathcal{B}'$. For each
$b'\in \mathcal{B}'$ (or $b\in\mathcal{B}$) define a sequence
$(w^{(b)}_n)_n = (w^{(b')}_n)_n$ of positive numbers in the following 
way: \ If $b'=\emptyset$ then
$w^{(b')}_n = 1/n$ for all $n$. Let $b'=(b'_n)\in\mathcal{B}'$ be a non-empty
finite or infinite sequence of positive integers. Set $w^{(b')}_n = 1/b'_1$ for
$n=1,\ldots, b'_1$. Assume that $w^{(b')}_i$ has been defined for
$i=1,\ldots,  b'_1+\ldots+b'_n$ and $w^{(b')}_{b'_1+\ldots+b'_n}=1/K$
for some $K\in \mathbb{N}$.
If $b'=(b'_1,\ldots, b'_n)$ then $w^{(b')}_i = 1/(K+i)$ for $i>b'_1
+\ldots+b'_n$. Let $b'_1 = (b'_1,\ldots, b'_n, b'_{n+1},\ldots)$. If
$b'_{n+1}>K$ then $w^{(b')}_i=1/b'_{n+1}$ for
$i=b'_1 +\ldots+ b'_n+1,\ldots, b'_1 +\ldots+b'_{n+1}$. If 
$b'_{n+1}\le K$ then $w^{(b')}_i =1/(K+i)$ for
$i=b'_1 +\ldots+ b'_n+1,\ldots, b'_1 +\ldots+ b'_{n+1}$.
\\
Thus for example if $b'=(2,3,2)$ then
$w^{(b')}=(\frac12,\frac12,\frac13,\frac13,\frac13,\frac14,
\frac14,\frac14,\frac14,\frac15,\frac16,\frac17,\ldots).$

Note that for every
$b'\in \mathcal{B}'$
$$w^{(b')}_n\searrow 0\quad \text{and}\quad \sum_n w^{(b')}_n = \infty.$$
Also note that the sequences $(w^{(b')}_n)_n$ for $b'\in\mathcal{B}'$ satisfy
the following properties
\begin{equation}\label{eq22}
\left\{\begin{matrix}
\text{for every sequence $(\lambda_n)$ with $\lim_n\lambda_n=\infty$
there
exists }
b'\in\mathcal{B}' \text{ such that}\hfill\\
\hfill \lim\limits_{k\to\infty} \frac{\sum\limits^k_{n=1} 
w^{(b')}_n}{\lambda_k}
=0\hfill
\end{matrix}\right.
\end{equation}
and
\begin{equation}\label{eq23}
\left\{\begin{matrix}
\text{for all $b'=(b'_n)\in\mathcal{B}'$ and for all $m\in\mathbb{N}$ there
exists }
M\in\mathbb{N} \text{ such that}\hfill\\
  \text{if $\beta'=(\beta'_n)\in \mathcal{B}'$
with $\beta'_n= b'_n$ for }
n=1,\ldots, M, \text{ then}\hfill\\
\hfill w^{(\beta')}_m \le w^{(b')}_m.\hfill
\end{matrix}\right.
\end{equation}

In order to see (\ref{eq22}) choose $n_1 < n_2<\ldots$ such that
$$\inf_{i\ge n_k} \lambda_i \ge (k+1)^2 \quad \text{and}\quad n_{k+1} > n_1
+\ldots+n_k\quad \text{for}\quad k\in\mathbb{N}.$$
Let $b' = (n_1,n_2,\ldots)\in \mathcal{B}'$. Then
$$\frac{\sum\limits^k_{n=1} w^{(b')}_n}{\lambda_k} \le \frac1k \quad
\text{for}\quad k\ge n_1$$
which proves (\ref{eq22}). Property (\ref{eq23}) is obvious.

Let $n\in \mathbb{N}$ and let $s_1,\ldots, s_n$, $t_1,\ldots, t_n\in T$ such
that
\begin{itemize}
\item[$\bullet$] $s_1 \prec s_2 \prec\ldots\prec s_n$
\item[$\bullet$] $s_i \preceq t_i$ for $i=1,\ldots, n$
\item[$\bullet$] $\{t\in T\colon \ s_i \preceq  t \preceq t_i\} \cap \{t\in
T\colon \ s_j \preceq t \preceq t_j\} = \emptyset$ for $i\ne j$.
\item[$\bullet$] If $s_i\prec t_i$ for some $i$ then $\{t\in T\colon \ s_i
\prec t \preceq t_i\} \cap \{t\in T\colon \ s_1 \preceq t \preceq s_n\} =
\emptyset$.
\end{itemize}
Then the set
$$\mathcal{S} = \{t\in T\colon \ s_1 \preceq t \preceq s_n\} \cup
\bigcup^n_{i=1} \{t\in T\colon \ s_i \preceq t \preceq t_i\}$$
is called a {\it bush\/} and it is denoted by bush$(s_1,t_1,s_2,t_2, \ldots,
s_n,t_n)$. For such a bush $\mathcal{S}$ we write $|\mathcal{S}|=n$.

We shall define $Y$ to be the completion of $c_{00}(T)$, the linear
space of finitely supported real valued functions on $T$, under a
certain norm given below. For $\alpha \in T$ we let $e_\alpha \in
c_{00}(T)$ be defined by $e_\alpha (\gamma) = \delta_{ \alpha ,
\gamma}$. $(e_\alpha)_{\alpha \in T}$ will be a basis for $Y$ in the
lexicographic order, i.e. $\{ e_{(0,0)}, e_{(1,0)}, e_{(1,1)}, e_{(2,0)} ,
\ldots \}$. $(e_\alpha ^*)_{\alpha \in T}$ will denote the
biorthogonal functionals of the basis. First, if $k\in\mathbb{N}\cup
\{0\}$  and $\mathcal{S} =
\text{bush}(s_1,t_1,\ldots, s_n,t_n)$ then we define a seminorm
$\|\cdot\|_{k,\mathcal{S}}$ on $c_{00}$ as follows
$$\|x\|_{k,\mathcal{S}} = \sup \left\{\sum^n_{j=1} w^{(b)}_{k+j}
|e^*_{t_j}(x)|\colon  b\in\mathcal{B},\ s_n\in b\right\}.$$
For $x\in c_{00}$ we define
$$\|x\| = \|x\|_\infty \vee \sup (\|x\|_{0,\mathcal{S}_1} +
\|x\|_{|\mathcal{S}_1|, \mathcal{S}_2} +\ldots+ \|x\|_{|\mathcal{S}_1| +
|\mathcal{S}_2| + \ldots+ |\mathcal{S}_{n-1}|, \mathcal{S}_n})$$
where the supremum is taken over $n\in\mathbb{N}$ and mutually disjoint bushes
$\mathcal{S}_1,\mathcal{S}_2,\ldots, \mathcal{S}_n$. Then $(Y,\|\cdot\|)$ is
the completion of $(c_{00}(T),\|\cdot\|)$. It is obvious that
$(e_t)_{t\in T}$ is
a normalized basis for $Y$. The Banach space $Y$ satisfies the following

\begin{Prop}\label{pro3.4}
\begin{itemize}
\item[a)] For every sequence $(\lambda_n)_n$ with $\lim_n \lambda_n = \infty$
there exists a subsequence $(x_i)$ of $(e_\alpha)$ having a
spreading model $(\tilde x_i)$ such that
$$\lim_n \frac{\left\|\sum\limits^n_{i=1} \tilde x_i\right\|}{\lambda_n} = 0.$$
\item[b)] For every normalized block basis $(x_i)$ of $(e_\alpha)$ having a
spreading model $(\tilde x_i)$ we have that
$$\lim_n \left\|\sum^n_{i=1} \tilde x_i\right\| = \infty.$$
\end{itemize}
\end{Prop}

\begin{proof}
a) \ First choose $b' = (b'_n)\in \mathcal{B}'$ so that (\ref{eq22})
is satisfied. Recall that $b'$ corresponds to a branch $b = \{u_1
\prec u_2 \prec \ldots\}\in\mathcal{B}$ of $T$. Let $(x_i)=(e_{u_i})$
and let $(w_n^{(b')})=(w_n^{(b)})= (w_n)$. We claim that for every
$m\in \mathbb{N}$
\begin{equation}\label{eq24}
\limsup_{\stackrel{\scriptstyle i_1 <i_2<\ldots<i_m}{\scriptstyle
i_1\to \infty}} \left\|\sum^m_{k=1} x_{i_k}\right\| \le \sum^m_{k=1} w_k+1
\end{equation}
which will finish the proof of a) by passing to a subsequence of
$(x_i)$ with a spreading model. By (\ref{eq23}) there exists $M\in
\mathbb{N}$ such that if $\beta \in \mathcal{B}$ and $u_M \in \beta$
then $w_n^\beta \leq w_n$ for $n=1, \ldots , m$. Let $i_1 < i_2 <
\ldots < i_m$ with $M < i_1$, and set
$$x = \sum^m_{k=1} x_{i_k}.$$ Let $\mathcal{S}_i$, $i=1,\ldots, N$ be
disjoint bushes with $\mathcal{S}_i = \text{bush}(s_{n_{i-1}+1},
t_{n_{i-1}+1},\ldots, s_{n_i},t_{n_i})$ where $0 = n_0 < n_1 < \ldots
< n_N$ and assume that $\{t_{n_{i-1}+1},\ldots, t_{n_i}\} \cap
\{u_{i_1},\ldots, u_{i_m} \}\ne \emptyset$ for all $i=1,\ldots, N$. By
the definition of the bush there exists at most one $i\in\{1,\ldots,
N\}$ such that $s_{n_i} \preceq u_{i_1}$. Without loss of generality
assume that $s_{n_1} \preceq u_{i_1}$. Then $\{t_1,\ldots,
t_{n_1-1}\}\cap \{u_{i_1},\ldots, u_{i_m}\} = \emptyset$ and
\begin{equation}\label{eq25}
\|x\|_{k,\mathcal{S}_1} \le 1
\end{equation}
for all $k \in \mathbb{N} \cup \{ 0 \}$. Also for $i=2,\ldots, N$ we
have that $u_{i_1}\preceq s_{n_i}$. Thus for $i=2,\ldots, N$ and $k\in
\mathbb{N}$ there exists $\beta\in \mathcal{B}$ with $s_{n_i}\in
\beta$ such that
\begin{equation}\label{eq26}
\|x\|_{k,\mathcal{S}_i} = \sum^{K_i}_{j=1} w^{(\beta)}_{k+j} \le
\sum^{K_i}_{j=1} w^{(b)}_{k+j}
\end{equation}
where $K_i = \bigl|
\{t_{n_{i-1}+1},\ldots, t_{n_i}\}\cap \{u_{i_1},\ldots,
u_{i_m}\}\bigr|
\wedge (m-k)$. Combining (\ref{eq25}) and (\ref{eq26}) we obtain
$$\|x\| \le 1 + \sum^m_{k=1} w^{(b)}_k$$
which implies (\ref{eq24}).

\noindent b) \ Let a normalized block basis  $(x_i)$ of $(e_\alpha)$ have a
spreading model  $(\tilde x_i)$. We consider two
cases:
\noindent {\it Case 1:} \ Assume that $\|x_i\|_\infty \not\to 0$.
By passing to a subsequence we can assume without loss of generality that
there exists a sequence $(t_n)$ in $T$ and $\varepsilon>0$ such that
$$|e^*_{t_n}(x_n)| >\varepsilon \quad \text{for all}\quad n.$$ Since
$(x_n)$ is a block sequence, we have that $|t_n|\to \infty$ (where for
$t\in T$ we set $|t|= \bigl|\{s\in T\colon \ (0,0) \preceq s \preceq
t\}\bigr|$). Thus it is easy to show that there exists a subsequence
$(t_{i_n})$ of $(t_n)$ and there exists $(s_n)\subset T$ with $s_1
\prec s_2 \prec \ldots$ and $s_n \preceq t_{i_n}$ for all $n$. We can
further assume without loss of generality that
$$
\{t\in T\colon \ s_n \preceq t \preceq t_{i_n}\} \cap \{t\in T\colon \ s_m
\preceq t \preceq t_{i_m}\} = \emptyset\quad \text{for}\quad n\ne m
$$
  and for $n\in \mathbb{N}$ if
$s_n \prec t_{i_n}$ then
$$\{t\in T\colon \ s_n \prec t \preceq t_{i_n} \} \cap b = \emptyset$$
where $b\in\mathcal{B}$ is the branch that contains all the $s_i$'s. Thus for
all $n\in \mathbb{N}$ and $j_1 < j_2 <\ldots<j_n$,
\begin{align*}
\left\|\sum^N_{k=1} x_{i_{j_k}}\right\| &\ge \left\|\sum^n_{k=1}
x_{i_{j_k}} \right\|_{0, \text{ bush}(s_{j_1}, t_{i_{j_1}},\ldots,
s_{j_n}, t_{i_{j_n}})}\\ &\ge \sum^n_{k=1}
w^{(b)}_k|e^*_{t_{i_{j_k}}}(x_{i_{j_k}})| \ge \sum^n_{k=1}
w^{(b)}_k\varepsilon \mathop{\longrightarrow}_{n\to \infty} \infty.
\end{align*}

\noindent {\it Case 2:} \ Assume that $\|x_i\|_\infty\to 0$.
In this case we use the following result whose proof is postponed.

\begin{Lem}\label{lem3.5}
For all $K\in\mathbb{N}$ there exists $\varepsilon>0$ such that for
all $x\in Y$ with $\|x\|=1$ and $\|x\|_\infty<\varepsilon$ there exist
$n\in\mathbb{N}$ and disjoint bushes, $\mathcal{S}_1,
\mathcal{S}_2,\ldots, \mathcal{S}_n$ with
$$.5 < \|x\|_{K,\mathcal{S}_1} + \|x\|_{K+|\mathcal{S}_1|,\mathcal{S}_2}
+\ldots+ \|x\|_{K+|\mathcal{S}_1|+\ldots+ |\mathcal{S}_{n-1}|,
\mathcal{S}_n}.$$
\end{Lem}

In this case, by Lemma \ref{lem3.5} we obtain a subsequence $(y_i)$ of
$(x_i)$ such that for all $i \in \mathbb{N}$ there exist disjoint
bushes $\mathcal{S}^{(i)}_1,\ldots, \mathcal{S}^{(i)}_{\ell_i}$ such
that
$$
.5 < \|y_i\|_{K_{i-1},\mathcal{S}^{(i)}_1} +\ldots+ \|y_i\|_{K_{i-1}
+|\mathcal{S}^{(i)}_1| +\ldots+ |\mathcal{S}^{(i)}_{\ell_i-1}|,
\mathcal{S}^{(i)}_{\ell_i}}
$$
where $K_0=0$ and $K_{i-1}=|\mathcal{S}^{(1)}_1| + \ldots +
|\mathcal{S}^{(i-1)}_{\ell_{i-1}}|$ for $i>1$. Thus for all $n\in
\mathbb{N}$  and $j_1 < j_2 <\ldots< j_n$ we have
\begin{align*}
\left\|\sum^n_{k=1} y_{j_k}\right\| &\ge \sum^n_{k=1} \left(
\|y_{j_k}\|_{K_{j_k-1}, \mathcal{S}^{(j_k)}_1} +\ldots+
\|y_{j_k}\|_{K_{j_k-1} + |\mathcal{S}^{(j_k)}_1| +\ldots+
|\mathcal{S}^{(j_k)}_{\ell_{j_k}-1}|,
\mathcal{S}^{(j_k)}_{\ell_{j_k}}} \right) \ge .5 \, n
\mathop{\longrightarrow}_{n\to\infty} \infty
\end{align*}
which finishes the proof of Proposition \ref{pro3.4}.
\end{proof}

It only remains to give the

\begin{proof}[Proof of Lemma \ref{lem3.5}]
Choose $\varepsilon>0$ that satisfies $K\varepsilon<.4$. Then for
$x\in Y$ with $\|x\|=1$ and $\|x\|_\infty<\varepsilon$, let
$\widetilde{\mathcal{S}}_1,\ldots, \widetilde{\mathcal{S}}_m$ be
disjoint bushes with
\begin{equation}\label{eq27}
.9 < \|x\|_{0,\widetilde{\mathcal{S}}_1} +
\|x\|_{|\widetilde{\mathcal{S}}_1|, \widetilde{\mathcal{S}}_2}
+\ldots+ \|x\|_{|\widetilde{\mathcal{S}}_1| +\ldots+
|\widetilde{\mathcal{S}}_{m-1}|, \widetilde{\mathcal{S}}_m}.
\end{equation}
Let $i_0$ be the smallest integer with $|\widetilde{\mathcal{S}}_1|
+\ldots+ |\widetilde{\mathcal{S}}_{i_0}| \ge K$. Let
$\widetilde{\mathcal{S}}_{i_0,1}$, $\widetilde{\mathcal{S}}_{i_0,2}$
be two disjoint bushes, subsets of $\widetilde{\mathcal{S}}_{i_0}$,
with $|\widetilde{\mathcal{S}}_{i_0,1}| +
|\widetilde{\mathcal{S}}_{i_0,2}|=|\widetilde{\mathcal{S}}_{i_0}|$ and
$|\widetilde{\mathcal{S}}_1| +\ldots+
|\widetilde{\mathcal{S}}_{i_0-1}| + |\widetilde{\mathcal{S}}_{i_0,1}|
= K$ (if $i_0=1$ then $|\widetilde{\mathcal{S}}_1| +\ldots+
|\widetilde{\mathcal{S}}_{i_0-1}| = 0$). Note that
\begin{equation}\label{eq28}
\|x\|_{0,\widetilde{\mathcal{S}}_1} +\ldots+
\|x\|_{|\widetilde{\mathcal{S}}_1| +\ldots+
|\widetilde{\mathcal{S}}_{i_0-1}|, \widetilde{\mathcal{S}}_{i_0,1}}
\le
(|\widetilde{\mathcal{S}}_1| +\ldots+
|\widetilde{\mathcal{S}}_{i_0-1}| + |\widetilde{\mathcal{S}}_{i_0,1}|)
\|x\|_\infty
  \le K\varepsilon \nonumber
  < .4.\nonumber
\end{equation}
By (\ref{eq27}) and (\ref{eq28}) we have that the choice $\mathcal{S}_1 =
\widetilde{\mathcal{S}}_{i_0,2}$, $\mathcal{S}_2 =
\widetilde{\mathcal{S}}_{i_0+1},\ldots, \mathcal{S}_n =
\widetilde{\mathcal{S}}_m$ $(n=m-i_0+1)$ satisfies the statement of
Lemma~\ref{lem3.5}.
\end{proof}

For all $n \in \N$ it is easy to construct a space $X$ for which the
cardinality $|SP(X)| = |SP_\omega (X)| = n$. Indeed, $X =
\bigl(\sum_{i=1}^n \ell_{p_i}\bigr)_2$ suffices, where the $p_i$'s are
distinct elements of $(1, \infty)$. Also if $2 < p_1 < p_2 < \ldots$
then it is not hard to show that
$|SP_\omega\bigl(\bigl(\sum_{i=1}^\infty \ell_{p_i}\bigr)_2\bigr)| =
\omega$.  In this case one obtains an infinite decreasing chain of
spreading models.

But we do
   not know very much what happens hereditarily. Let us mention
  some questions (among many) concerning the "heriditary structure of 
spreading models".

\begin{Qst} \label{qst3.3}
Does there exist a Banach space such that in every infinite
dimensional subspace there exist normalized basic sequences
having spreading models equivalent to the unit vector bases of
$\ell_1$ and $\ell_2$? If such space exists, must it contain more
(perhaps uncountably many) mutually non-equivalent spreading models?
More generally, does there exist $X$ so that for all subspaces $Y$ of
$X$ and $1 \le p < \infty$, the unit vector basis of $\ell_p$ (and of
$c_0$) is equivalent to a spreading model of $Y$? Is the space
constructed in \cite{OS} or \cite{OS2} such a space?
\end{Qst}

In order to answer Question~\ref{qst3.3}, the answer to the following
question may be useful:

\begin{Qst}  \label{qst3.4}
Can we always isomorphically (or isometrically) stabilize
the set of spreading models by passing to appropriate subspaces?
i.e. for every Banach space $X$ does there exists a subspace $Y$ such
that for every normalized basic sequence $(y_i)$ in $Y$ having
spreading model $(\tilde{y}_i)$ and for every further subspace $Z$ of
$Y$, there exists a normalized basic sequence $(z_i)$ in $Z$ having
spreading model $(\tilde{z}_i)$ such that $(\tilde{z}_i)$ is
equivalent (respectively, isometric) to $(\tilde{y}_i)$? Is the space
$X$ constructed in section \ref{sec2} a counterexample?
\end{Qst}

\begin{Qst} \label{qst3.5}
  Let $n\in\N$. Does there exist a Banach space so that every subspace 
has exactly
  $n$ (isomorphically or isometrically) different spreading models?
  Does there exist a Banach space so that every subspace has exactly
  countably many(isomorphically or isometrically) different spreading models?
\end{Qst}
Many problems exist concerning the structure of the partially
ordered set $SP_\omega(X)$ (in the sense of Definition \ref{D3.1}). 
We state a few of these.  If
$SP_\omega(X)$ is countable, then $SP_\omega(X)$ must have an upper
bound by Proposition~\ref{mainlemma}.

\begin{Qst} \label{qst3.6} What are the realizable isomorphic 
structures of the partially
     ordered set $(SP_\omega(X), \le)$? In particular, for every finite
     partially ordered set $(P, \le)$ such that any two elements admit
a least  upper bound, does there
     exist $X$ such that $SP_\omega(X)$ is isomorphic to $(P, \le)$.
\end{Qst}
  \begin{Qst} \label{qst3.7}    Suppose $SP_\omega(X)$ is finite (or 
even countable). What can
     be said about $X$? Must some spreading model be equivalent to the
     unit vector basis in $c_0$ or $\ell_p$ ($1 \le p < \infty$)?  We
     address the case $|SP_\omega(X)| =1$ in Section~4.
\end{Qst}

\section{Spaces with a unique spreading model}\label{sec4}

\begin{Qst}\label{qst4.1}
Let $X$ be an infinite dimensional Banach space so that $|SP(X)|=1$.
Must the unique spreading model of $X$ be equivalent to the unit
vector basis of $\ell_p$ for some $1\le p <\infty$, or $c_0$?
\end{Qst}

One could also raise similar questions by restricting to either those
spreading models generated by normalized weakly null basic sequences
or, in the case that $X$ has a basis, to those generated by normalized
block bases.

We give some partial answers to these questions using our techniques above.

\begin{Prop}\label{pro4.2}
Let $X$ be a Banach space an infinite dimensional Banach space so that
all  spreading
models of normalized basic sequences in $X$ are equivalent.
\begin{itemize}
\item[a)]  If all the spreading models are uniformly equivalent, i.e.
  if there exists $D\in \mathbb{R}$ so that the spreading
models of all normalized basic sequences in $X$ are $D$-equivalent,
then all spreading models of $X$ are equivalent to the unit vector
basis of $\ell_p$ for some $1\le p < \infty$ or $c_0$.
\item[b)] Let $(z_i)$ be a normalized basic sequence which dominates a
(hence every) spreading model of $X$. Then there exists $C<\infty$ so
that $(z_i)$ $C$-dominates any spreading model of a normalized basic
sequence $(x_i)$ in $X$.
\item[c)] If $p$ belongs to the Krivine set of the spreading model $(\tilde
x_i)$ of some normalized basic sequence $(x_i)$ of $X$ then
$(\tilde x_i)$ dominates the unit vector basis of $\ell_p$.
\item[d)] If 1 belongs to the Krivine set of some spreading model in
$X$ then all spreading models are equivalent to the unit vector basis
of $\ell_1$.
\end{itemize}
\end{Prop}

\begin{proof} If $X$ is not reflexive then there exists a normalized basic
sequence $(x_n)$ in $X$ which dominates the summing basis \cite{J}. By
\cite{R} $(x_n)$ has a subsequence $(x_{n_k})$ which is either
equivalent to the unit vector basis of $\ell_1$ or it is
weak-Cauchy. In the later case $(x_{n_{2k+1}} - x_{n_{2k}})_k$ is
weakly null and thus by passing to a subsequence we can assume that it
has an unconditional spreading model which dominates the summing basis
and hence must be equivalent to the unit vector basis of
$\ell_1$. Therefore in either case there exists a spreading model in
$X$ equivalent to the unit vector basis of $\ell_1$, and it is easy to
see that a)--d) hold. Thus for the proof of a)--d) we may assume that
$X$ is reflexive.

a) Let $(\tilde x_i)$ be a spreading model of $X$ and let $p$ in the
Krivine set of $(\tilde x_i)$. By Remark~\ref{rem1.1} for every $n\in
\mathbb{N}$ there exists a spreading model $(\tilde x^{(n)}_i)_i$ of
$X$ such that $(\tilde x^{(n)}_i)^n_{i=1}$ is 2-equivalent to the unit
vector basis of $\ell^n_p$. Also $(\tilde x_i)^n_{i=1}$ is
$D$-equivalent to $(\tilde x^{(n)}_i)^n_{i=1}$ thus $2D$-equivalent to
the unit vector basis of $\ell^n_p$.

b) Let $(z_i)$ be a normalized basic sequence which dominates all
spreading models of $X$.  Assume that the statement is false.  Then
for every $n \in \N$ there exists a normalized weakly null basic
sequence $(x^n_i)$ in $X$, having spreading model $(\tilde{x}^n_i)$,
and there exist scalars $(a^n_i)_i$, such that $\| \sum_i a^n_i
\tilde{x}^n_i \| = 2^{2n}$ and $\| \sum_i a^n_i z_i \| =1$. By
Proposition~\ref{mainlemma} there exists a seminormalized weakly null
sequence $(y_i)$ in $X$, having spreading model $(\tilde{y}_i)$ such
that $(\tilde{y}_i)$ $2^n$-dominates $(\tilde{x}^n_i)$ for all
$n$. Thus $\| \sum_i a^n_i \tilde{y}_i \| \ge 2^{-n} \| \sum_i a^n_i
\tilde{x}^n_i \| = 2^n$. Hence $(z_i)$ does not dominate
$(\tilde{y}_i)$, which is a contradiction.

c) This follows from b) and Remark \ref{rem1.1}.

d) This follows  from c).
\end{proof}

\begin{Rmk} \label{rmk4.3}
If $X$ a basis $(e_i)$ and the hypothesis of Proposition \ref{pro4.2}
is changed to ``all spreading models of normalized block bases are
equivalent'' then one obtains a similar theorem, while the conclusions
are restricted to spreading models generated by normalized block
bases. The ``$X$ is not reflexive'' part of the proof is replaced by
``$(e_i)$ is not shrinking''. If the hypothesis is changed to ``all
spreading models generated by normalized weakly null basic sequences
are equivalent'' then one has two cases: Either $X$ is a Schur space,
hence $X$ is hereditarily $\ell_1$ \cite{R}, or $X$ does admit such a
spreading model. And the proposition holds in the latter case with the
obvious modifications.
\end{Rmk}

A problem closely related to 4.1 has been considered by
Ferenczi, Pelczar and Rosendal in \cite{FPR}:
Suppose that $X$ has a basis $(e_i)$ for which every normalized block basis
has a subsequence equivalent to $(e_i)$.
Must $(e_i)$ be equivalent to the unit vector basis of $c_0$ or some $\ell_p$?
Using proposition~3.1, the authors obtain results analogous to those in
this paper.

Many additional questions remain about the structure of the spreading
models of a Banach space $X$.

\section{Existence of non-trivial operators on subspaces of
  certain Banach spaces} \label{sec5}

In this section we give sufficient conditions on a Banach space $X$
for the existence of a subspace $Y$ of $X$ and an operator $T:Y \to X$
which is not a compact perturbation of a multiple of the inclusion map.
This property is related to the long standing open problem of whether
there exists a Banach space (of infinite dimension)
on which every operator is a
compact perturbation of a multiple of the identity. Notice that if a
Banach space $X$ contains an unconditional basic sequence then there
exists a subspace $Y$ of $X$ and an operator $T: Y \to Y$ such that
$T^n$ is not a compact perturbation of a multiple of the identity for
all positive integers $n$ ($Y$
can be taken to be the closed linear span of the unconditional basic
sequence, and $T$ to be the projection on an infinite
subsequence). W.T. Gowers \cite{G} proved that there exists a subspace
$Y$ of the Gowers-Maurey space $GM$
(as defined in \cite{GM}), and an operator $T: Y \to GM$ which is not
a compact perturbation of a multiple of the inclusion. In \cite{AS} it is
shown that there exists an operator on $GM$ which is not a compact
perturbation of a
multiple of the identity. It is also known that the asymptotic
$\ell_1$ Hereditarily
Indecomposable (H.I.) space constructed by S.A. Argyros and I. Deliyanni
\cite{AD} admits such operator (unpublished work of S.A. Argyros
and R. Wagner). Certain asymptotic $\ell_1$ H.I. spaces constructed
by I. Gasparis \cite{Ga1} admit such operators as well
\cite{Ga2}. Our approach generalizes  the idea of \cite{G}.

\begin{Thm} \label{Main3}
   Let $X$ be a Banach space. Assume that there exists a normalized
   weakly null basic sequence $(x_i) $ in $X$ having spreading model
   $(\tilde{x}_i)$ which is not equivalent to the unit vector basis of
   $\ell_1$ yet $1$ belongs to the Krivine set of $(\tilde{x}_i)$. Then
   there exists a subspace $W$ of $X$ and a continuous linear operator
   $T: W \to W$ such that $p(T)$ is not a compact perturbation of a
   multiple of the identity operator on $W$, for every non-constant
   polynomial $p$.
\end{Thm}

The proof uses a convenient auxiliary  notation.
%%\begin{Def}
Let ${\mathcal F} \subseteq [ \N ]^{< \infty }$ be a family of finite
subsets of  positive integers.  For $(a_i) \in c_{00}$ we set
$$
\| (a_i) \|_{\ell_1({\mathcal F})}= \sup \{ \sum_{i \in F} |a_i|: F
\in {\mathcal F} \} .
$$
%%\end{Def}

\medskip

\begin{proof}[Proof of Theorem \ref{Main3}]
The main part of the proof is the following

\noindent
{\bf Claim 1:} For every
$\ell \in \N \cup \{ 0 \}$ there exists $(w^\ell_i)_i$ a
seminormalized sequence in $X$, an increasing sequence
$(M^{(\ell +1)}_i)_i$ of positive integers and a sequence
$(\delta^{(\ell +1)}_i)_i$ of positive numbers with
$\sum_i \delta^{(\ell +1)}_i < \infty$, such that
$w^0_1, w^1_1, w^0_2, w^2_1, w^1_2, w^0_3, \ldots$ is a basic sequence
in  $X$,
and  for every
$(a^{(\ell)}_j)_{ \ell \in \N \cup \{ 0 \} , j \in \N }
\in c_{00}((\N \cup \{ 0 \} )\times \N)$
we have
\begin{equation} \label{34}
\max_{1 \leq \ell < \infty } \| (a^{(\ell)}_j)_j \|_\ell \leq
\| \sum_{\ell =0}^\infty \sum_{j=1}^\infty a^{(\ell)}_j w^\ell_j \|
\leq \sum_{ \ell=0}^\infty \| (a^{(\ell)}_j)_j \|_{\ell +1}.
\end{equation}
Here for $\ell \in \N $ and $(a_j)_j \in c_{00}$ we put
$$
\| (a_j )_j \|_\ell := \sup_{i \in \N} \delta^{(\ell)}_i \| (a_j)_j
\|_{\ell_1({\mathcal G}^\ell_i)}
$$
where for $\ell , i \in \N$ we set
${\mathcal G}^\ell _i =
\{ F \subset \N : |F| \leq M^{(\ell)} _i \} $.

\medskip

For $\ell=0$ we consider the left-hand side expression in
(\ref{34}) to be equal to 0, and we adopt the same convention
throughout  the rest of this section.

\medskip

Once Claim 1 is established, let $\widetilde{W}=\mbox{span}\{
w^\ell_j : \ell \in \N \cup \{ 0 \} ,\ j \in \N \}$ and define
$T: \widetilde{W} \to \widetilde{W}$ by
$$
T (w^0_j)= 0 \mbox{ and } T(w^{\ell + 1}_j )= \frac{1}{2^{\ell +1}}w^\ell
_j \mbox{ for all } \ell \in \N \cup \{ 0 \} \mbox{ and } j \in \N .
$$
Since $(w^n_i)_{n \in \N \cup \{ 0 \}, i \in \N}$ is a basic sequence
in  $X$, $T$ is well defined. Let
$(a^{(\ell)}_j)_{ \ell \in \N \cup \{ 0 \} , j \in \N }
\in c_{00}((\N\cup \{ 0 \} ) \times \N)$
and
$x = \sum_{\ell =0}^\infty \sum_{j=1}^\infty
a^{(\ell)}_j w^\ell_j \in\widetilde{W}$.
We have
\begin{align*}
\| T x \| = & \| \sum_{\ell =1}^\infty \sum_{j=1}^\infty a^{(\ell)}_j
\frac{1}{2^\ell} w^{\ell -1}_j \|
\le \sum_{\ell =1}^\infty \frac{1}{2^\ell}
\| (a^{(\ell)}_j)_j \|_\ell \mbox{ (by (\ref{34}))} \\
\leq & \max_{1 \le \ell < \infty} \| (a^\ell_j)_j \|_\ell \\
\le & \| x \| \mbox{ (by (\ref{34}))}.
\end{align*}
Thus if $W$ denotes the closure of $\widetilde{W}$ then $T$ extends
to a bounded operator on $W$.

Let $p(t)=a_n t^n + a_{n-1}t^{n-1} + \cdots + a_1 t +a_0$ be a
non-constant  polynomial.
We show that $p(T)$ is not a compact perturbation of a
multiple of the identity operator $I$ on $W$.
Indeed, for any $i \in \{ 1,2, \ldots , n\}$ and $j \in \N$ we have
$$
T^iw^n_j = T^{i-1} \frac{1}{2^n}w^{n-1}_j =
T^{i-2} \frac{1}{2^n} \frac{1}{2^{n-1}}w^{n-2}_j = \cdots =
\prod_{k=n-i+1}^n \frac{1}{2^k} w^{n-i}_j.
$$
Thus for every scalar
$\lambda$ and $j \in \N$ we have:
$$
(p(T)- \lambda I)w^n_j= (\sum_{i=1}^n a_i T^i - \lambda I)w^n_j=
\sum_{i=1}^n a_i \left( \prod_{k=n-i+1}^n \frac{1}{2^k} \right)
w^{n-i}_j + (a_0- \lambda)w^n_j.
$$
Since $w^0_1, w^1_1, w^0_2, w^2_1, w^1_2, w^0_3, \ldots$ is a basic sequence
in  $X$, there exist $j_1< j_2< \cdots $ in $\N$  such
that
$(\sum_{i=1}^n a_i ( \prod_{k=n-i+1}^n \frac{1}{2^k} )
w^{n-i}_j + (a_0- \lambda)w^n_j)_j$
is a seminormalized block sequence of
$w^0_1, w^1_1, w^0_2, w^2_1, w^1_2, w^0_3, \ldots$ which proves that
$p(T) - \lambda I$ is not a compact operator.

\medskip
Claim 1 follows from

\noindent
{\bf Claim 2} There exists a subspace $Y$ of $X$ with a basis and for
every $\ell \in \N \cup \{ 0 \}$ there exists a seminormalized weakly
null basic sequence $(u^\ell_i)_i$ in $Y$, an increasing sequence
$(M^{(\ell +1)}_i)_i$ of positive integers, and a sequence
$(\delta^{(\ell +1)}_i)_i$ of positive numbers with $\sum_i
\delta^{(\ell +1)}_i < \infty$, such that the vectors $u^n_i$ for $n
\in \N \cup \{ 0 \}$ and $i \in \N$ are disjointly supported with
respect to the basis in $Y$, and for every $\ell \in \N \cup \{ 0 \}$
and $(a^{(m)}_j)_{m \in \{ 0,1, \ldots , \ell \} , j \in \N} \in
c_{00}(\{ 0, 1, \ldots , \ell \} \times \N)$
%and $y \in Y$
%such that the vectors $y$ and
%$\sum_{m=0}^\ell \sum_{j=1}^\infty a^{(m)}_j u^m_j$ are disjointly
%supported with respect to the basis of $Y$,
we have that
\begin{equation} \label{37}
  2 \max_{1 \leq m \leq \ell} \| (a^{(m)}_j)_j \|_m \leq
\| \sum_{m=0}^\ell \sum_{j=1}^\infty a^{(m)}_j u^m_j \|
\leq (1/2) \sum_{m=0}^\ell \| (a^{(m)}_j)_j \|_{m+1}.
\end{equation}

Once Claim 2 is established, passing for every $n \in \N$ to a
subsequence of $(u^n_i)_i$ and making small perturbations if
necessary, we get $(w^n_i)_i$ such that $w^0_1, w^1_1, w^0_2, w^2_1,
w^1_2, w^0_3, \ldots$ forms a block basis in $Y$ and for all
$(a^{(\ell)}_j)_{ \ell \in \N \cup \{ 0 \} , j \in \N } \in c_{00}((\N
\cup \{ 0 \})\times \N)$ we have
\begin{equation} \label{39}
\frac{1}{2} \| \sum_{\ell =0}^\infty \sum_{j=1}^\infty
a^{(\ell)}_j u^\ell_j \|
\leq \| \sum_{\ell =0}^\infty \sum_{j=1}^\infty a^{(\ell)}_j w^\ell_j \|
\leq 2 \| \sum_{\ell =0}^\infty \sum_{j=1}^\infty
a^{(\ell)}_j u^\ell_j \|.
\end{equation}
Obviously (\ref{37}) and (\ref{39}) imply (\ref{34}) and thus
Claim 1 follows.

\medskip

Now we prove Claim 2. We construct the space $Y$ and inductively on
$\ell \in \N \cup \{ 0 \}$ we construct the sequences $(u^\ell_i)_i$,
$(M^{(\ell+1)}_i)_i$ and $(\delta^{(\ell +1)}_i)$ which satisfy
(\ref{37}). The upper and lower estimates are based on the following
two lemmas of independent interest, whose proofs we postpone until the
end of the section.

In their formulation, given an increasing sequence $M_1<M_2< \cdots$
of integers, we let ${\mathcal G}_n = \{ G \subset \N :  |G| \leq M_n
\}$ for $n \in \N$.

\medskip

\begin{Lem} \label{nonell1}
   Let $X$ be a Banach space and $(x_i)_i$ be a normalized weakly null
   basic sequence in $X$ which has a spreading model $(\tilde{x}_i)$
   not equivalent to the unit vector basis of $\ell_1$. Then for every
   $(\delta_n)_{n \ge 2} \subset (0, 1)$ there exists a subsequence
   $(x_{m_i})$ of $(x_i)$, an increasing sequence $M_1<M_2< \cdots$ of
   integers, and $\delta_1>0$ such that for all $(a_i)_i \in c_{00}$ we
   have
\begin{equation} \label{upperbound}
\| \sum a_i x_{m_i} \| \leq \sup_{n \in \N} \delta_n \| (a_i)_i
\|_{\ell_1({\mathcal G}_n)}.
\end{equation}
\end{Lem}

\medskip

\begin{Lem} \label{Kriv}
   Let $X$ be a Banach space and $(z_i)_i$ be a normalized weakly null
   basic sequence in $X$ which has spreading model $(\tilde{z}_i)_i$
   such that $1$ belongs to the Krivine set of $(\tilde{z}_i)_i$.
   There exists a subsequence $(z_i')_i$ of $(z_i)_i$ with the
   following property.  Given any infinite subset $J \subseteq \N$, any
   subsequence $(M_n)_n$ of $\N$, and $(\delta_n)_n \subset (0,\infty)$
   with $\sum_{n=1}^\infty \delta_n < \infty$, there exists a
   seminormalized weakly null basic sequence $(y_i)_i$ in the span of
   $(z_j')_{j \in J}$ of disjointly supported vectors with respect to
   $(z_j')_{j \in J}$ such that for all $(a_i) \in c_{00}$ and all $y $
   in the span of $(z_j')_{j \not\in J}$ we have
%%%%
%\begin{equation} \label{K}
%\sup_{ n \in \N} \delta_n \| (a_i)
%\|_{\ell_1({\mathcal G}_n)} \leq \|  \sum a_i y_{k_i} \| .
%\end{equation}
%Moreover,  for any normalized weakly null basic sequence $(z^0_i)$ in $X$
%there exists a subsequence $(y^0_i)_i$ of $(z^0_i)$ and there exist
%normalizedvectors $(y^n_i)_{n, i \in \N}$ in $X$ such that
%$y^0_1, y^1_1, y^0_2, y^2_1, y^1_2, y^0_3, y^3_1, \ldots$ is a
%basic sequence in $X$ spanning a space $Y$, such that
%for every $(\delta_n)_n$ and $({\mathcal G}_n)_n$
%as above, and for every
%$K=\{ s_1 < s_2 < \cdots \} \in [ \N ]$ there exists a
%seminormalized sequence $(y_i)$ in $Y$ such that
%$y_i \in \mbox{span}(y^n_{s_i})_{n \in \N}$ for all $i$, satisfying:
%for all  $(a_i) \in c_{00}$, $k_1<k_2< \cdots $ in $\N$,  and for
%all vectors  $y \in Y$ such that $y$ and $\sum a_i y_{k_i}$ are
%disjointly supported with respect to
%$(y^n_i)_{n \in \N \cup \{ 0 \} , i \in \N }$,  we have that
%%%%
\begin{equation} \label{K'}
\sup_{ n \in \N} \delta_n \| (a_i)
\|_{\ell_1({\mathcal G}_n)} \leq \|  y+ \sum a_i y_{i} \|.
\end{equation}
Furthermore, if $(\tilde{z}_i)$ is not equivalent to the unit vector
basis of $\ell_1$ then no spreading model of $(y_i)$ is equivalent to
the unit vector basis of $\ell_1$.
\end{Lem}

We now return  to the proof of Claim 2.

\medskip

Since $1$ belongs to the Krivine set of $(\tilde{x}_i)$, pass to a
subsequence of $({x}_i)_i$ constructed in Lemma~\ref{Kriv}.  Without
confusion denote this subsequence again by $({x}_i)_i$, and let $Y$
be its span.  Let $K_0, K_1, K_2, \ldots \in [ \N ]$ be disjoint
subsequences of positive integers. For all $\ell \in \N \cup \{0\}$ we
will construct $u^\ell_i \in \mbox{span}\{ x_j: j \in K_\ell \}$.

\medskip

\noindent
{\em Construction of  $(u^0_i)_i$, $(\delta^{(1)}_i)_{i \in \N}$
   and $(M^{(1)}_i)_{i \in \N}$:}
Let  $(\delta^{(1)}_i)_{i \ge 2} \subset (0, 1)$ such that
$\sum_{i \geq 2} \delta^{(1)}_i < \infty$. Since $(\tilde{x}_i)_i$ is
a spreading model of $(x_j)_{j \in {K_0}}$
which is not
isomorphic to the unit vector basis of $\ell_1$ we may
apply Lemma \ref{nonell1} to obtain a
subsequence $(x_{m_i})$ of $(x_j)_{j \in {K_0}}$,
an increasing sequence
$(M^{(1)}_i)_{i \in \N}$
of positive integers, and $\delta^{(1)}_1>0 $ such that for
all $(a_i) \in c_{00}$ we have
\begin{equation} \label{40}
\| \sum a_i x_{m_i} \| \leq \frac12 \| (a_i) \|_1: = \frac12 \sup_{n \in \N}
\delta^1_n \| (a_i) \|_{\ell_1({\mathcal G}^1_n)}
\end{equation}
where ${\mathcal G}^1_n = \{ G \subset \N : |G|\leq M^{(1)}_n \}$.
This is (\ref{37}) for $\ell =0$.

\medskip
\noindent
{\em The inductive step - Construction of $(u^\ell_i)_i$,
$(\delta^{(\ell+1)}_i)_i$ and $(M^{(\ell+1)}_i)$:} Assume that we have
constructed $(u^m_i)_i$, $(M^{(m+1)}_i)_i$ and $(\delta^{(m+1)}_i)$
for $m=0, 1, \ldots , \ell -1$ so that (\ref{37}) is satisfied when
$\ell$ is replaced by $\ell -1$.
Apply Lemma~\ref{Kriv} for $J = K_\ell$,
$(M_i)_i=(M^{(\ell)}_i)_i$ and $(\delta_i)_i=
(\delta^{(\ell)}_i)_i$ to obtain
a seminormalized weakly null basic
sequence $(u^\ell_i)_i$ in $\mbox{span}\{ x_j: j \in K_\ell \}$
satisfying: for all $(a_i) \in
c_{00}$,  and
$y \in \mbox{span}\{ x_j: j \not\in K_\ell \}$
we have that
\begin{equation} \label{42}
2 \sup_{ n \in \N} \delta^{(\ell)}_n \| (a_i)
\|_{\ell_1({\mathcal G}^\ell_n)} \leq \|  y+ \sum a_i u^\ell_i \|,
\end{equation}
where ${\mathcal G}^\ell_n = \{ G \subset \N : |G|\leq M^{(\ell)}_n
\}$ for $n \in \N$.
By passing to a subsequence of $(u^\ell_i)_i$ and
relabeling we can assume that $(u^\ell_i)_i $ has a spreading model
$(\tilde{u}^\ell_i)$. By the ``furthermore'' part of Lemma
\ref{Kriv} we have that $(\tilde{u}^\ell_i)$ is not equivalent to the
unit vector basis of $\ell_1$. Let
$(\delta^{(\ell +1)}_i)_{i \ge 2} \subset (0,1)$ such that
$\sum_{i \ge 2} \delta^{(\ell +1)}_i < \infty$. Apply Lemma
\ref{nonell1} to obtain a subsequence of $(u^\ell_i)_i$ (which we
still call $(u^\ell_i)_i$) , an increasing sequence
$(M^{(\ell+1)}_i)_{i \in \N}$ of positive integers, and
$\delta^{(\ell +1)}_1>0$ such that for all $(a_i) \in c_{00}$ we have
\begin{equation} \label{43}
\| \sum a_i u^\ell_i \| \leq \frac12 \| (a_i ) \|_{\ell +1} :=
\frac12 \sup_{n \in \N} \delta^{(\ell +1)}_n
\| (a_i) \| _{\ell_1({\mathcal G}^{\ell+1}_n)}
\end{equation}
where ${\mathcal G}^{\ell +1}_n = \{ G \subset \N :  |G| \leq M^{(\ell
+1)}_n \}$ for $n \in \N$.

We now show that (\ref{37}) is satisfied.
Let
$(a^{(m)}_j)_{m \in \{ 0,1, \ldots , \ell \} , j \in \N} \in c_{00}(\{
0, 1, \ldots , \ell \} \times \N)$
By the triangle inequality and the
inductive hypothesis we can estimate $\| \sum_{m=0}^\ell
\sum_{j=1}^\infty a^{(m)}_j u^m_j \|$ as follows:
\begin{align} \label{44}
2 \max_{1 \leq m \le \ell -1} \| (a^{(m)}_j)_j \|_m \leq &
\|  \sum_{j=1}^\infty a^{(\ell)}_j u^\ell_j  +
\sum_{m=0}^{\ell -1} \sum_{j=1}^\infty a^{(m)}_j u^m_j \| \\
\leq & \frac12 \sum_{m=0}^{\ell -1} \| (a^{(m)}_j)_j \|_{m+1} +
\| \sum_{j=1}^\infty a^{(\ell)}_j u^\ell_j \|  \nonumber \\
%%\le &\frac12 \sum_{m=0}^{\ell -1} \| (a^{(m)}_j)_j \|_{m+1} +
%% \| \sum_{j=1}^\infty a^{(\ell)}_j u^\ell_j \| \nonumber \\
\le &\frac12 \sum_{m=0}^{\ell -1} \| (a^{(m)}_j)_j \|_{m+1} +
  \frac12 \| (a_i^{(\ell)} ) \|_{\ell +1} \mbox{ (by (\ref{43}))}
\nonumber \\
= & \frac12 \sum_{m=0}^{\ell} \| (a^{(m)}_j)_j \|_{m+1}. \nonumber
\end{align}
By (\ref{42}) we also  have
%$\| y + \sum_{m=0}^\ell \sum_{j=1}^\infty a^{(m)}_j u^m_j \|$ from below
%as follows:
\begin{equation} \label{49}
2 \sup_{ n \in \N} \delta^{(\ell)}_n \| (a^{(\ell)}_i)
\|_{\ell_1({\mathcal G}^\ell_n)} \leq
\|  \sum_{m=0}^{\ell -1} \sum_{j =1}^\infty a^{(m)}_j u^m_j
+ \sum_{j=1}^\infty a^{(\ell)}_j u^\ell_j \|.
\end{equation}
Obviously (\ref{37}) follows immediately from  (\ref{44}) and (\ref{49}).
\end{proof}

If we are interested only in the construction of an operator on a
subspace which is not a compact perturbation of a multiple of the
inclusion map, then the spreading model assumptions of
Theorem~\ref{Main3} can be significantly relaxed and the argument
would be essentially simpler.

\begin{Thm} \label{Main2}
Let $X$ be a Banach space. Assume that
there exist normalized weakly null basic sequences $(x_i)$, $(z_i)$
in $X$ such that $(x_i)$ has spreading
model $(\tilde{x}_i)$ which is not equivalent to the unit vector basis
of $\ell_1$, and $(z_i)$ has spreading model $(\tilde{z}_i)$ such that
$1$ belongs to the Krivine set of $(\tilde{z}_i)$. Then there exists a
subspace $Y$ of $X$ and an
operator $T:Y \to X$ which is not a compact
perturbation of a multiple of the inclusion map.
\end{Thm}

\begin{proof}[Sketch of  proof]
Let $(\delta_n)_{ n \ge 2} \subset (0, 1)$ such that
$\sum_{n=2}^\infty \delta_n < \infty$.
Using  Lemma~\ref{nonell1}
we obtain  a subsequence
$(x_{m_i})_i$ of $(x_i)$, an increasing sequence
$M_1<M_2<\cdots$ of integers and $\delta_1 >0$
such that for all $(a_i) \in c_{00}$ we
have
%%\begin{equation} \label{35}
$$
\| \sum a_i x_{m_i} \| \le \sup_{n \in \N} \delta_n
\| (a_i)_i \|_{\ell_1({\mathcal G}_n)},
$$
%%\end{equation}
where ${\mathcal G}_n = \{ G \subset \N :  |G| \leq M_n \}$
for $n \in\N$. Then by Lemma~\ref{Kriv}
we  obtain a seminormalized weakly null basic
sequence $(y_i)$ in the span of
$(z_i)$ such that for all
$(a_i) \in c_{00}$ and $k_1<k_2< \ldots$ in $\N$,
%%\begin{equation} \label{lowerest}
$$
\| \sum a_i y_{k_i} \| \geq \sup_{ n \in \N} \delta_n \| (a_i)
\|_{\ell_1({\mathcal G}_n)}.
$$
%%\end{equation}
Thus for every $(a_i) \in c_{00}$ we have $ \| \sum a_i x_{m_i} \| \le
\| \sum a_i y_i \|$, and passing to subsequences if necessary we may
also assume that $x_{m_1}, y_1, x_{m_2}, y_2, \ldots$ is a
(seminormalized weakly null)  basic sequence.
Thus the operator $T$ defined on $[y_i: i \in \N]$, the closed linear
span of $(y_i)$, by $T(y_i)=x_{m_i}$ for all $i$, is a continuous
operator. Also for any scalar $\lambda$ the
operator $T- \lambda I$ (where $I$ denotes the inclusion operator
from $[y_i: i \in \N ]$ to $X$) is non-compact, since
$(T-\lambda I)(y_i)= x_{m_i}- \lambda y_i$ which is a seminormalized
weakly null sequence.
\end{proof}

We now give the proofs of Lemmas \ref{nonell1} and \ref{Kriv}.

\begin{proof}[Proof of Lemma~\ref{nonell1}]
Since $(x_i)$ is a normalized weakly null sequence, by a slight 
generalization of a theorem
  of J.~Elton
\cite{E} (see \cite{DKK}),
  by passing  to a subsequence and relabeling we can
assume the following. There exists $C>0$  such that for every $(a_i) \in
c_{00}$ with $|a_i|\le 2$, for all $i\in\N$, and for every $n\in\N$ 
and $F \subseteq \{ i: 2^{-n-1}\le |a_i| < 2^{-n} \}$ we
have that
\begin{equation} \label{g1}
\| \sum_{i \in F}a_i x_i \| \le C \| \sum a_i x_i \| .
\end{equation}
  We can also assume that for every $n \in \N$, if $F
\subseteq \{ n , n+1, \ldots \}$, $|F| \leq n 2^n$ and $(a_i) \in
c_{00}$ then
\begin{equation} \label{g3}
\| \sum_{i \in F}  a_i x_i \| \le 2 \| \sum_{i \in F} a_i
\tilde{x}_i \| .
\end{equation}
Since $(\tilde{x}_i)$ is not equivalent to the unit vector basis
of $\ell_1$ there exists $\varepsilon_n \searrow 0$ such that for
all $m$ we have
\begin{equation} \label{g4}
\| \sum_{i=1}^m \tilde{x}_i \| \le \varepsilon_m m .
\end{equation}
Let $n_0 \in \N$ such that $\sum_{n \geq n_0} n 2^{-n} < 1/8$,
${\mathcal F}_0 = \{ F \subseteq \N : | F | \leq
\sum_{i=1}^{n_0}i2^i \}$, and for $n \in \N$ let ${\mathcal F}_n =
\{ F \subseteq \N : | F| \leq n 2^n \} $. Now let $x = \sum a_i
x_i$ with $\| x \| =1$. Let $F_0 = \{ i : 2\ge |a_i | >2^{-1} \}$
and for $n \in \N$ let $F_n = \{ i : 2^{-(n+1)} \leq |a_i | <
2^{-n} \}$. We have
\begin{align} \label{g5}
\| x \| &\leq \sum_{n=0}^\infty\limits \| \sum_{i \in F_n }
\limits a_i x_i \|  \nonumber \\ &= \sum_{ \{ n : |F_n| > n2^n \}
}\limits \| \sum_{i \in F_n } \limits a_i x_i \| + \sum_{ \{ n
\geq n_0 : | F_n | \leq n2^n \} } \limits \| \sum_{i \in F_n }
\limits a_i x_i \| + \sum_{ \{ n < n_0 : | F_n | \leq n 2^n \} }
\limits \| \sum_{i \in F_n } \limits a_i x_i \| .
\end{align}
Notice that by (\ref{g1}) if $| F_n | \geq n2^n$ then
\begin{equation} \label{g6}
\| \sum_{i \in F_n } a_i x_i \| \leq C = 2 C n^{-1} 2^{-(n+1)}n
2^n \leq 2 C n^{-1} \| (a_i) \|_{\ell_1({\mathcal F}_n)} .
\end{equation}
Let $n \geq n_0$ with $|F_n| \leq n2^n$. Let $x^* \in X^*$ with
$\| x^* \| =1$ such that $\| \sum_{i \in F_n } a_i x_i \| = x^*(
\sum_{i \in F_n} a_i x_i )$. Let $F_{n,1}= \{ i \in F_n: a_i >0,
x^*(x_i)>0 \}$, $F_{n,2}= \{ i \in F_n: a_i >0, x^*(x_i)<0 \}$,
$F_{n,3}= \{ i \in F_n : a_i <0, x^*(x_i) >0 \}$, $F_{n,4} = \{ i
\in F_n : a_i <0, x^*(x_i)<0 \}$. There exists $j_n \in \{ 1,
\ldots 4 \}$ and $y^* \in \{ x^*, -x^* \}$ such that if we set
$F_n'= F_{n,j_n}$ we have
\begin{equation*}
\| \sum_{i \in F_n} a_i x_i \| \leq 4 | \sum_{i \in F_n'} a_i
x^*(x_i) | \leq 4 \sum_{i \in F_n'} | a_i ||x^*(x_i)| \leq 4 \
2^{-n} y^*( \sum_{i \in F_n'} x_i ) \leq 4 \ 2^{-n} \| \sum_{i \in
F_n'} x_i \| .
\end{equation*}
Thus if $|F_n'| \leq n$ then $\| \sum_{i \in F_n } a_i x_i \| \leq
4n2^{-n}$. If $|F_n' | >n $ then since $|F_n'| \leq |F_n| \leq
n2^n$, by (\ref{g3}) and (\ref{g4}) we have
\begin{align*}
\| \sum_{i \in F_n} a_i x_i \| & \leq 4 \ 2^{-n} \| \sum_{i \in
F_n'} x_i \| \leq 8 \ 2^{-n} \| \sum_{i \in F_n'} \tilde{x}_i \|
\\ &\leq 8 \ 2^{-n} \varepsilon_n |F_n| \leq 16 \varepsilon_n
\sum_{i \in F_n'} | a_i | \leq 16 \varepsilon_n \| (a_i)
\|_{\ell_1({\mathcal F}_n)} .
\end{align*}
Hence
\begin{align} \label{g7}
\sum_{ \{ n \geq n_0: |F_n| \leq n2^n \} } \| \sum_{i \in F_n} a_i
x_i \| & \leq \sum_{  n \geq n_0  } 4n 2^{-n} +
\sum_{  n \geq n_0 } 16 \varepsilon_n  \| (a_i) \|_{\ell_1({\mathcal
F}_n)} \nonumber \\ & \leq 4 \sum_{n \geq n_0} n2^{-n} + \sum_{ \{ n
\geq n_0  } 16 \varepsilon_n \| (a_i)
\|_{\ell_1({\mathcal F}_n)} .
\end{align}
Finally,
\begin{equation} \label{g8}
\sum_{ \{ n <n_0 : | F_n| \leq n2^n \} } \| \sum_{i \in F_n } a_i
x_i \| \leq \sum_{ \{ n<n_0 : |F_n| \leq n2^n \} } \sum_{i \in
F_n} | a_i| \leq \| (a_i) \|_{\ell_1({\mathcal F}_0)} .
\end{equation}
By (\ref{g5})-(\ref{g8}) and the choice of $n_0$ we have
\begin{align*}
\| x \| & \leq \sum_{n=1}^\infty 2C n^{-1} \| (a_i)
\|_{\ell_1({\mathcal F}_n)} + 4 \sum_{n \geq n_0} n 2^{-n} +
\sum_{  n>n_0} 16 \varepsilon_n \| (a_i)
\|_{\ell_1({\mathcal F}_n)} + \| (a_i) \|_{\ell_1({\mathcal F}_0)}
\\ & < \frac12 + \sum_{n=0}^\infty \varepsilon_n' \| (a_i)
\|_{\ell_1({\mathcal F}_n)}
\end{align*}
where $\varepsilon_0'=1$ and $\varepsilon_n'= \max (2Cn^{-1}, 16\varepsilon_n
)$ for $n \in \N$. Since $\| x \| =1$ we obtain
\begin{equation*}
\| x \| \leq \sum_{ n=0 }^ \infty 2 \varepsilon_n' \| (a_i)
\|_{\ell_1({\mathcal F}_n)} .
\end{equation*}
Note that if $k_1<k_2< \cdots$ in $\N$, $M_n= \max \{ |F|: F \in
{\mathcal F}_n \}$ (for $n \in \N \cup \{ 0 \} $), ${\mathcal
G}_0= \{ F \subseteq \N : |F| \leq M_0 + \cdots + M_{k_1} \}$,
${\mathcal G}_1 = \{ F \subseteq \N: |F| \leq M_{k_1+1} + \cdots +
M_{k_2} \} , \ldots$, $\tilde{\varepsilon}_0= 2 \max \{
\varepsilon_n': n \in \{ 0, \ldots , k_1 \} \}$,
$\tilde{\varepsilon}_1= 2 \max \{ \varepsilon_n': n \in \{ k_1 +1,
\ldots , k_2 \} \} , \ldots$, then
\begin{equation*}
\| x \| \leq \sum_{n=0}^\infty \tilde{\varepsilon}_n \| (a_i)
\|_{\ell_1({\mathcal G}_n)} .
\end{equation*}
Thus, if we choose $k_1<k_2< \cdots$ such that $2^{n+1}
\tilde{\varepsilon}_n \leq \delta_n$ for $n \in \N$ we shall have
\begin{equation*}
\sum_{n=0}^\infty \tilde{\varepsilon}_n \| (a_i)
\|_{\ell_1({\mathcal G}_n)} \leq \sup_{ n \in \N \cup \{ 0 \} }
2^{n+1} \tilde{\varepsilon}_n \| (a_i) \|_{\ell_1({\mathcal G}_n)}
\leq \sup_{n \in \N \cup \{ 0 \} } \delta_n \| (a_i) \|_{\ell_1
({\mathcal G}_n)} ,
\end{equation*}
which finishes the proof of Lemma \ref{nonell1}.
\end{proof}

\medskip

\begin{proof}[Proof of Lemma \ref{Kriv}]
%   Fix $\epsilon >0$.
   Since $1$ belongs to the Krivine set of $(\tilde z_i)$, pick for
   every $n \in \N$ vectors $\tilde w_j^n \in \mbox{span} (\tilde
   z_i)$, for $j=1, 2, \ldots$ such that for any subset $E \subseteq
   \N$ with $|E| = n$, vectors $(\tilde w_j^n)_{j \in E}$ are
   $2$ equivalent to the unit vector basis in $\ell_1^n$.
   Moreover, $\tilde w_1^n, \tilde w_2^n, \ldots$ are successive and
   have the same distribution with respect to $(\tilde z_i)$. In
   particular, denote the common length of their support by $K_n$.

   Fix $J$, $(M_n)$ and $(\delta_n)$ as in the assumptions.  Using
   Schreier unconditionality theorem (\cite{MR}, also \cite{BL} or
   \cite{O}) pass to a subsequence $(z_i')$ of $(z_i)$ such that for
   any finite subset $F \subseteq \N$ such that $|F| \le n K_n $ and $n
   < \min F$, for some $n \in \N$, we have
\begin{equation}
   \label{schreier}
   \|\sum_{i \in F} a_i z_i'\| \le 3\|\sum_{i}
   a_i z_i'\|,
\end{equation}
for any scalars $(a_i)$.

For any $n \in \N$, let $(w_j^n)$ be equidistributed vectors in the
span of $(z_i')_{i \in J}$ corresponding to $(w_j^n)$, and supported
after $z_n'$, by the definition of a spreading model (\ref{eq0}). In
particular we may require that for any subset $E \subseteq \N$ with
$|E| = n$, the sequence $(\tilde w_j^n)_{j \in E}$ is $3$ equivalent
to the unit vector basis in $\ell_1^n$.  We may additionally chose the
$w_j^n$'s so that $w_1^1, w_1^2, w_2^1, w_1^3, \ldots$ form a block
basis with respect to $(z_i')_{i \in J}$.

For j=1, 2, \ldots, set
$$
y_j = \sum_n \delta_n  w_j^{M_n}.
$$
We clearly have $\max_n \delta_n \le \|y_j\| \le \sum_n \delta_n$,
for all $j$.

To prove (\ref{K'}), pick $(a_i) \in c_{00}$ and a vector $y$
supported outside of $J$.  Fix $n \in \N$ and $G \subseteq \N $ with
$|G| \le M_n$.  Noting that the supports of $w_j^n$'s with respect to
$(z_i')$ have cardinality $K_n$, by Schreier unconditionality
(\ref{schreier}) we can isolate $w_j^n$'s from the expression for
$y_j$ to get
$$
%%\begin{eqnarray*}
\|y + \sum_j a_j y_j \| \ge
(1/3) \|\sum_{j \in G} a_j \delta_n  w_j^n \|
\ge  (1/9)\delta_n \sum_{j \in G} |a_j|.
%%\end{eqnarray*}
$$
Taking into account the definition of ${\mathcal G}_n$ and of the norm
$ \| \cdot\|_{\ell_1({\mathcal G}_n)}$  this completes the proof.
\end{proof}

It is proved in \cite{AS} that the spreading model of the unit
vector basis of the
Gowers-Maurey space $GM$ as defined in \cite{GM}, is
isometric to the unit vector basis
of Schlumprecht space $S$ as defined in \cite{S2}. Thus
Theorem \ref{Main3} immediately gives the following:

\begin{Cor} \label{Gowers}
There exists a subspace $Y$ of $GM$ and an operator $T: Y \to Y$
such that $p(T)$ is not a compact perturbation of a multiple of the
identity, for every non-constant polynomial $p$.
\end{Cor}

%%%% OLD

\vspace{.5in}
\scriptsize{
\noindent
Department of Mathematics, University of South Carolina, Columbia, SC
29208.
giorgis@math.sc.edu

\noindent
Department of Mathematics, University of Texas at Austin,  Austin,
TX 78712.
odell@math.utexas.edu

\noindent
Department of Mathematics, Texas A\&M University, College Station,
TX 77843.
schlump@math.tamu.edu

\noindent
Department of Mathematical Sciences, University of Alberta, Edmonton,
Alberta, T6G 2G1, Canada. \\ nicole@ellpspace.math.ualberta.ca
}

\end{document}